\documentclass[11pt]{article}
\usepackage[margin=3.2cm]{geometry}

\usepackage{latexsym,amssymb}
\usepackage{amsmath}
\usepackage{color}
\usepackage{mathscinet}

\newtheorem{Theorem}{\bf Theorem}[section]
\newtheorem{Lemma}{\bf Lemma}[section]
\newtheorem{Proposition}{\bf Proposition}[section]
\newtheorem{Corollary}{\bf Corollary}[section]
\newtheorem{Remark}{\bf Remark}[section]
\newtheorem{Example}{\bf Example}[section]
\newtheorem{Definition}{\bf Definition}[section]

\newenvironment{theorem}{\begin{Theorem}$\!\!\!$}{\end{Theorem}}
\newenvironment{lemma}{\begin{Lemma}$\!\!\!$}{\end{Lemma}}
\newenvironment{proposition}{\begin{Proposition}$\!\!\!$}{\end{Proposition}}
\newenvironment{corollary}{\begin{Corollary}$\!\!\!$}{\end{Corollary}}
\newenvironment{remark}{\begin{Remark}$\!\!\!$}{\end{Remark}}

\newenvironment{definition}{\begin{Definition}$\!\!\!$}{\end{Definition}}

\numberwithin{equation}{section}

\begin{document}

\title{Existence of solutions for an inhomogeneous\\ fractional semilinear heat equation}
\author{
\qquad\\
Kotaro Hisa, Kazuhiro Ishige and Jin Takahashi
}
\date{}
\maketitle
\begin{abstract}
We obtain necessary conditions and sufficient conditions 
on the existence of solutions to the Cauchy problem for a fractional semilinear heat equation 
with an inhomogeneous term. 
We identify the strongest spatial singularity 
of the inhomogeneous term for the solvability of the Cauchy problem.
\end{abstract}
\vspace{25pt}
\noindent Addresses:

\smallskip
\noindent K.~H.:  Mathematical Institute, Tohoku University,\\
\qquad\,\,\,\,\, 6-3 Aoba, Aramaki, Aoba-ku, Sendai 980-8578, Japan. \\
\noindent 
E-mail: {\tt kotaro.hisa.s5@dc.tohoku.ac.jp}\\

\smallskip
\noindent 
K.~I.:  Graduate School of Mathematical Sciences, The University of Tokyo,\\
\qquad\,\,\, 3-8-1 Komaba, Meguro-ku, Tokyo 153-8914, Japan. \\
\noindent 
E-mail: {\tt ishige@ms.u-tokyo.ac.jp}\\

\smallskip
\noindent J. T.:  Department of Mathematical and Computing Science, Tokyo Institute of Technology,\\
\qquad\,\,\, 2-12-1 Ookayama, Meguro-ku, Tokyo 152-8552, Japan. \\
\noindent 
E-mail: {\tt takahashi@c.titech.ac.jp}\\
\vspace{20pt}
\newline
\noindent
{\it MSC:} 35A01; 35K58; 35R11
\vspace{3pt}
\newline
\noindent
{\it Keywords:} Semilinear heat equation; Fractional Laplacian; Inhomogeneous term
\vspace{3pt}
\newline
\newpage

%

\section{Introduction}
This paper is concerned with the Cauchy problem for a fractional semilinear heat equation 
with an inhomogeneous term 
\begin{equation}\label{eq:1.1}
	\left\{
	\begin{aligned}
	&\partial_t u+(-\Delta)^{\frac{\theta}{2}}u = u^p + \mu,\quad && x\in{\bf R}^N,\,\,\, t>0, \\
	&u(0) = 0 && \mbox{in}\quad {\bf R}^N,
	\end{aligned}
	\right.
\end{equation}
where $\partial_t:=\partial/\partial t$, $N\geq1$, $0<\theta\leq 2$, $p>1$ 
and $\mu$ is a nonnegative Radon measure in ${\bf R}^N$ 
or a nonnegative measurable function in ${\bf R}^N$. 
Here $(-\Delta)^{\theta/2}$ denotes the fractional power of the Laplace operator $-\Delta$ in ${\bf R}^N$. 
In this paper we study necessary conditions and sufficient conditions 
on the inhomogeneous term $\mu$ for the existence of nonnegative solutions to problem~\eqref{eq:1.1} 
and identify the strongest singularity of $\mu$ for the solvability of problem~\eqref{eq:1.1}.  
Our identification is new even for $\theta=2$. 

Before considering problem~\eqref{eq:1.1}, 
we recall some results on the Cauchy problem 
\begin{equation}\label{eq:1.2}
	\left\{
	\begin{aligned}
	&\partial_t u+(-\Delta)^{\frac{\theta}{2}}u = u^p, \quad && x\in{\bf R}^N,\,\,\, t>0, \\
	&u(0) = \nu && \mbox{in}\quad {\bf R}^N,
	\end{aligned}
	\right.
\end{equation}
where $N\ge1$, $0<\theta\leq 2$, $p>1$ and $\nu$ is a nonnegative Radon measure in ${\bf R}^N$ 
or a nonnegative measurable function in ${\bf R}^N$. 
In \cite{HI18} the first and the second authors 
of this paper 
studied necessary conditions and sufficient conditions on the initial data 
for the solvability of problem~\eqref{eq:1.2} and obtained the following property. 
\begin{itemize}
  \item[(a)] 
  Assume that there exists a nonnegative solution 
  to problem~\eqref{eq:1.2} in ${\bf R}^N\times[0,T)$ for some $T>0$. 
  Then there exists $c=c(N,\theta,p)>0$ such that 
  \begin{equation}
  \label{eq:1.3}
  \sup_{x\in{\bf R}^N}\nu(B(x,\sigma))\le c\sigma^{N-\frac{\theta}{p-1}},
  \quad
  0<\sigma\le T^{\frac{1}{\theta}}.
  \end{equation}
  Furthermore, 
  if $p=p_\theta:=1+\theta/N$, then there exists $c'=c'(N,\theta)>0$ such that 
  $$
  \sup_{x\in{\bf R}^N}\nu(B(x,\sigma))
  \le c'\left[\log\biggr(e+\frac{T^{\frac{1}{\theta}}}{\sigma}\biggr)\right]^{-\frac{N}{\theta}},
  \quad
  0<\sigma\le T^{\frac{1}{\theta}}.
  $$
  \end{itemize}
 Here $B(x,\sigma):=\{y\in{\bf R}^N\,:\,|x-y|<\sigma\}$ for $x\in{\bf R}^N$ and $\sigma>0$. 
 In the case of $1<p<p_\theta$, 
 since the function $\sigma\mapsto \sigma^{N-\theta/(p-1)}$ is decreasing for $\sigma>0$, 
 \eqref{eq:1.3} is equivalent to 
$$
\sup_{x\in{\bf R}^N}\nu(B(x,T^{\frac{1}{\theta}}))\le cT^{\frac{N}{\theta}-\frac{1}{p-1}}. 
$$
 Property~(a) implies 
  \begin{itemize}
  \item[(b)]
  There exists $c_1=c_1(N,\theta,p)>0$ such that, 
  if $\nu$ is a nonnegative measurable function in ${\bf R}^N$ satisfying 
  $$
  \begin{array}{ll}
  \nu(x)\ge c_1|x|^{-N}\displaystyle{\biggr[\log\left(e+\frac{1}{|x|}\right)\biggr]^{-\frac{N}{\theta}-1}}
  \quad & \mbox{if}\quad \displaystyle{p=p_\theta},\vspace{7pt}\\
  \nu(x)\ge c_1|x|^{-\frac{\theta}{p-1}}\quad & \mbox{if}\quad \displaystyle{p>p_\theta},
  \end{array}
  $$
  in a neighborhood of the origin, then problem~\eqref{eq:1.2} possesses no local-in-time solutions. 
  \end{itemize}
Furthermore, they obtained the following properties. 
  \begin{itemize}
  \item[(c)] Let $1<p<p_\theta$. Then there exists $c_2=c_2(N,\theta,p)>0$ such that, if
 $$
 \sup_{x\in{\bf R}^N}\nu(B(x,T^{\frac{1}{\theta}}))\le c_2T^{\frac{N}{\theta}-\frac{1}{p-1}}
 $$
 for some $T>0$, then problem~\eqref{eq:1.2} possesses a solution in ${\bf R}^N\times[0,T)$.
\item[(d)] Let $p\ge p_\theta$. Assume that 
$$
\begin{array}{ll}
0\le\nu(x)\le c|x|^{-N}\displaystyle{\biggr[\log\left(e+\frac{1}{|x|}\right)\biggr]^{-\frac{N}{\theta}-1}}+C_0
\quad & \mbox{if}\quad \displaystyle{p=p_\theta},\vspace{7pt}\\
0\le\nu(x)\le c|x|^{-\frac{\theta}{p-1}}+C_0\quad & \mbox{if}\quad \displaystyle{p>p_\theta},
\end{array}
$$
for some $c>0$ and $C_0\ge 0$. 
Then there exists $c_3=c_3(N,\theta,p)>0$ such that
problem~\eqref{eq:1.2} possesses a local-in-time solution if $c\le c_3$ 
and a global-in-time solution if $c\le c_3$, $C_0=0$ and $p>p_\theta$.
\end{itemize}
Assertions~(b) and (d) determine the strongest singularity of the initial data 
for the solvability of problem~\eqref{eq:1.2}. 
For related results, see e.g. \cite{BP85, HI19, IKO, IS16, RS13, Takahashi, We80}. 

On the other hand,
the existence of solutions to nonlinear parabolic equations with inhomogeneous terms 
has been studied in many papers, see e.g. \cite{BLZ00, BP85,  B, KT16, KT17, KK04, KQ02, L93, Zeng07, Zhang98_1, Zhang98_2, Zhang99} 
and references therein. However, there are no results concerning the identification 
of the strongest spatial singularity of the inhomogeneous term for the existence of solutions. 
In this paper, motivated by \cite{HI18}, 
we study necessary conditions and sufficient conditions on the inhomogeneous term~$\mu$  
for the existence of solutions to problem~\eqref{eq:1.1} and identify the strongest singularity 
of the inhomogeneous term $\mu$ for the solvability of \eqref{eq:1.1}. 
\vspace{3pt}

We formulate the definition of solutions to problem~\eqref{eq:1.1} and state our main results. 
\begin{definition}
\label{Definition:1.1}
Let $u$ be a nonnegative measurable function in ${\bf R}^N\times(0,T)$, where $0<T\le \infty$.  
We say that $u$ is a solution to problem~\eqref{eq:1.1} in ${\bf R}^N\times[0,T)$ if  
$u$ satisfies
\begin{equation*}
\int_0^T\int_{{\bf R}^N} u(-\partial_t \varphi + (-\Delta)^\frac{\theta}{2} \varphi )\,dx\,dt = 
\int_0^T\int_{{\bf R}^N}  u^p\varphi\,dx\,dt + \int_0^T \int_{{\bf R}^N}  \varphi \,d\mu(x)\,dt
\end{equation*}
for $\varphi\in C^\infty_0({\bf R}^N\times[0,T))$.
\end{definition}
The first theorem of this paper is concerned with necessary conditions 
on the inhomogeneous term $\mu$ for the solvability of problem~\eqref{eq:1.1}. 
Set 
$$
p_*:=\frac{N}{N-\theta}\quad\mbox{if}\quad 0<\theta<N
\quad\mbox{and}\quad
p_*:=\infty\quad\mbox{if}\quad \theta\ge N.
$$
\begin{theorem}
\label{Theorem:1.1}
Let $N\ge 1$, $0<\theta\le2$ and $p>1$. 
Let $u$ be a solution to problem~\eqref{eq:1.1} in ${\bf R}^N\times[0,T)$, 
where $0<T<\infty$. 
Then  there exists $\gamma=\gamma(N,\theta,p)>0$ 
such that
\begin{equation}
\label{eq:1.4}
\sup_{x\in{\bf R}^N}\mu(B(x,\sigma))\le \gamma\sigma^{N-\frac{\theta p}{p-1}}
\end{equation} 
for $0<\sigma\le T^{1/\theta}$. 
Furthermore, 
if $p=p_*$, then there exists $\gamma'=\gamma'(N,\theta)>0$ such that 
\begin{equation}
\label{eq:1.5}
\sup_{x\in{\bf R}^N}\mu(B(x,\sigma))
\le \gamma'\biggr[\log\biggr(e+\frac{T^{\frac{1}{\theta}}}{\sigma}\biggr)\biggr]^{-\frac{N}{\theta}+1}
\end{equation} 
for $0<\sigma\le T^{1/\theta}$.
\end{theorem}
If $1<p<p_*$, then the function $\sigma\mapsto \sigma^{N-\theta p/(p-1)}$ is decreasing for $\sigma>0$.
This means that 
\eqref{eq:1.4} is equivalent to 
$$
\sup_{x\in{\bf R}^N}\mu(B(x,T^\frac{1}{\theta}))\le \gamma\,T^{\frac{N}{\theta}-\frac{p}{p-1}}
$$
in the case of $1<p<p_*$. 
As corollaries of Theorem~\ref{Theorem:1.1}, we have
\begin{corollary}
\label{Corollary:1.1}
Let $N\ge 1$, $0<\theta\le2$ and $p\ge p_*$. 
Then there exists $\gamma=\gamma(N,\theta,p)>0$ such that, if 
a nonnegative measurable function $\mu$ in ${\bf R}^N$ satisfies 
$$
\mu(x)\ge 
\left\{\begin{aligned}
	&\gamma |x|^{-\frac{\theta p}{p-1}}&&\mbox{ if }\quad p>p_*,  \\
	&\gamma |x|^{-N} 
	\left[ \log \biggr( e + \frac{1}{|x|} \biggr) \right]^{-\frac{N}{\theta}}
	&&\mbox{ if }\quad p=p_*,
	\end{aligned}
	\right.
$$
in a neighborhood of the origin, then problem~\eqref{eq:1.1} possesses no local-in-time solutions. 
\end{corollary}
\begin{corollary}
\label{Corollary:1.2}
Let $N\ge 1$ and $0<\theta\le2$. 
\begin{itemize}
  \item[{\rm (1)}] 
  Let $1<p\le p_*$ and $\mu\not\equiv 0$ in ${\bf R}^N$. 
  Then problem~\eqref{eq:1.1} possesses no global-in-time solutions. 
  \item[{\rm (2)}] 
  Let $p>p_*$ and $\mu$ be a nonnegative measurable function in ${\bf R}^N$. 
  Then there exists $\gamma=\gamma(N,\theta,p)>0$ with the following property: 
  If there exists $R>0$ such that 
  $$
  \mu(x)\ge\gamma |x|^{-\frac{\theta p}{p-1}}
  $$
  for almost all $x\in{\bf R}^N\setminus B(0,R)$, then 
  problem~\eqref{eq:1.1} possesses no global-in-time solutions. 
  \end{itemize}
\end{corollary}
Next we state our results on sufficient conditions for the solvability. 
\begin{theorem}
\label{Theorem:1.2}
Let $N\geq1$, $0<\theta\leq 2$ and $1<p< p_*$. 
Then there exists $\gamma=\gamma(N,\theta,p)>0$ such that, 
if a nonnegative Radon measure $\mu$ in ${\bf R}^N$ satisfies 
$$
\sup_{x\in{\bf R}^N} \mu(B(x,\sigma)) 
\leq \gamma \sigma^{N-\frac{\theta p}{p-1}}
\quad \mbox{ for some }\sigma>0, 
$$
then problem \eqref{eq:1.1} possesses a solution in ${\bf R}^N\times[0,T)$ with $T=\sigma^\theta$. 
\end{theorem}
\begin{theorem}
\label{Theorem:1.3}
Let $N\geq1$, $0<\theta\leq 2$ and $p>p_*$. 
Let $1<r<\infty$ be such that
$$
r>r_*:=\frac{N(p-1)}{\theta p}.
$$
Then there exists $\gamma=\gamma(N,\theta,p,r)>0$ such that, 
if a nonnegative measurable function $\mu$ in ${\bf R}^N$ satisfies  
\begin{equation}
\label{eq:1.6}
\sup_{x\in{\bf R}^N}\|\mu\|_{L^r(B(x,\sigma))}\le \gamma \sigma^{\frac{N}{r}-\frac{\theta p}{p-1}}
\quad \mbox{ for some }\sigma>0, 
\end{equation}
then problem \eqref{eq:1.1} possesses a solution in ${\bf R}^N\times[0,T)$ with $T=\sigma^\theta$. 
\end{theorem}
%
\begin{theorem}
\label{Theorem:1.4}
Let $N\geq1$, $0<\theta\leq 2$ and $p\ge p_*$. 
Let $\mu$ be a nonnegative measurable function in ${\bf R}^N$ such that
\begin{equation}\label{eq:1.7}
	0\leq \mu(x)\leq 
	\left\{
	\begin{aligned}
	&\gamma |x|^{-\frac{\theta p}{p-1}} + C_0 &&\mbox{ if } \quad p>p_*,  \\
	&\gamma |x|^{-N} 
	\left[ \log \biggr( e + \frac{1}{|x|} \biggr) \right]^{-\frac{N}{\theta}} +C_0
	&&\mbox{ if } \quad p=p_*,
	\end{aligned}
	\right.
\end{equation}
for almost all $x\in {\bf R}^N$, where $\gamma>0$ and $C_0\ge 0$. 
Then there exists $\gamma_*=\gamma_*(N,\theta,p)>0$ such that 
problem~\eqref{eq:1.1} possesses a local-in-time solution if $\gamma\le\gamma_*$ 
and a global-in-time solution if $\gamma\le\gamma_*$, $C_0=0$ and $p>p_*$. 
\end{theorem}
By Theorems~\ref{Theorem:1.1}, \ref{Theorem:1.2} and \ref{Theorem:1.4} 
we can identify the strongest spatial singularity of $\mu$ for the solvability of problem~\eqref{eq:1.1}. 
Furthermore, 
by Theorems~\ref{Theorem:1.1} and \ref{Theorem:1.2} we easily obtain
\begin{corollary}
\label{Corollary:1.3}
Let $\delta$ be the Dirac delta function in ${\bf R}^N$. 
Then problem~\eqref{eq:1.1} possesses a local-in-time solution with $\mu=D\delta$ for some $D>0$ 
if and only if $1<p<p_*$. 
\end{corollary}
\begin{remark}
\label{Remark:1.1}
{\rm (i)} 
Corollary~{\rm\ref{Corollary:1.2}} {\rm (1)}
and Theorem~{\rm\ref{Theorem:1.4}}  
imply the following properties.
\begin{itemize}
  \item[{\rm (a)}] 
  If $1<p\le p_*$ and $\mu\not\equiv 0$, then problem~\eqref{eq:1.1} possesses no global-in-time solutions;
  \item[{\rm (b)}]
  If $p>p_*$, then problem~\eqref{eq:1.1} possesses a global-in-time solution for some $\mu\,(\not\equiv 0)$.
\end{itemize}
{\rm (ii)} 
In the case of $\theta=2$, 
assertions~{\rm (a)} and {\rm (b)} were first obtained by {\rm\cite{Zhang98_1}} 
and they have been extended to various nonlinear parabolic equations with inhomogeneous terms. 
See e.g. {\rm \cite{BLZ00, KK04, Zeng07, Zhang98_1, Zhang98_2, Zhang99}} and references therein. 
In the case of $0<\theta<2$, see {\rm \cite{KQ02}}.
\vspace{3pt}\newline
{\rm (iii)}  
Necessary conditions and sufficient conditions for the existence of solutions to the problem
\begin{equation*}
	\left\{
	\begin{aligned}
	&\partial_t u-\Delta u =|u|^{p-1}u + \delta\otimes\nu,\quad && x\in{\bf R}^N,\,\,\, t>0, \\
	&u(0) = 0 && \mbox{in}\quad {\bf R}^N,
	\end{aligned}
	\right.
\end{equation*}
were discussed in {\rm\cite{KT16, KT17}}, where $p>1$ and $\nu$ is a Radon measure in $[0,\infty)$. 
Corollary~{\rm\ref{Corollary:1.3}} with $\theta=2$ follows 
from {\rm\cite[{\it Theorem}~2.2]{KT16}} and {\rm\cite[{\it Theorem}~2.1]{KT17}}.
\end{remark}

We explain the idea of proving our theorems. 
Kartsatos and Kurta~\cite{KK04} obtained 
necessary conditions on the existence of global-in-time solutions to problem~\eqref{eq:1.1} with $\theta=2$.  
Except for the case of $0<\theta<2$ and $p=p_*$,
their arguments are available for the proof of Theorem~\ref{Theorem:1.1}. 
Indeed, the proof of Theorem~\ref{Theorem:1.1} except for such a case 
is given as a modification of the arguments in \cite{KK04}. 
In the case of $0<\theta<2$ and $p=p_*$, 
using a fractional Poisson equation,  
we modify  arguments in \cite{KK04} to prove Theorem~\ref{Theorem:1.1}. 
The regularity of solutions to the fractional Poisson equation plays an important role in the proof. 
On the other hand, the proofs of Theorems~\ref{Theorem:1.2} and \ref{Theorem:1.3} 
are based on the contraction mapping theorem in uniformly local Lebesgue spaces.  
Theorem~\ref{Theorem:1.4} is proved by the construction of supersolutions to problem~\eqref{eq:1.1}. 
This requires delicate estimates of volume potentials associated 
with the fundamental solution to the fractional heat equation.


The rest of this paper is organized as follows. 
In Section~2 we obtain necessary conditions 
for the solvability of problem~\eqref{eq:1.1} and prove Theorem~\ref{Theorem:1.1}. 
In Section~3 we apply the contraction mapping theorem in uniformly local Lebesgue spaces 
to prove Theorems~\ref{Theorem:1.2} and \ref{Theorem:1.3}. 
In Section~4 we prepare preliminary lemmas for the proof of Theorem~\ref{Theorem:1.4}. 
In Sections~5 and 6 we prove Theorem~\ref{Theorem:1.4} with $p>p_*$ and $p=p_*$, respectively. 
\section{Proof of Theorem~\ref{Theorem:1.1}}
We modify arguments in \cite{KK04} to prove Theorem~\ref{Theorem:1.1}.
We also prove Corollaries~\ref{Corollary:1.1} and \ref{Corollary:1.2}. 
In what follows, by the letter $C$
we denote generic positive constants and they may have different values also within the same line. 
For any set $E$ in ${\bf R}^N$, let $\chi_E$ be
 the characteristic function of $E$. 
\vspace{3pt}

We prepare the following lemma to prove \eqref{eq:1.4}. 
\begin{lemma}
\label{Lemma:2.1}
Let $u$ be a solution to problem~\eqref{eq:1.1} in ${\bf R}^N\times[0,T)$ for some $T>0$. 
Then, for any integer $s  > p/(p-1)$ and $z\in{\bf R}^N$, 
there exists $C=C(N,\theta,p,s)>0$ such that 
\begin{equation}
\label{eq:2.1}
\int_{{\bf R}^N} \zeta^s \, d\mu_z(x) 
\leq C T^{-\frac{p}{p-1}} \int_{{\bf R}^N} \zeta^s  \, dx 
+ C\int_{{\bf R}^N} |(-\Delta)^\frac{\theta}{2} \zeta|^\frac{p}{p-1}\zeta^{s-\frac{p}{p-1}} \, dx
\end{equation}
for $\zeta\in C^\infty_0({\bf R}^N)$ with $\zeta\geq0$. 
Here $\mu_z(A)=\mu(z+A)$ for Borel sets $A$ in ${\bf R}^N$. 
\end{lemma}
{\bf Proof.}
Let $u$ be a  solution to \eqref{eq:1.1} in ${\bf R}^N\times[0,T)$ with $0<T<\infty$. 
Let $z\in{\bf R}^N$ and set $u_z(x,t):=u(x-z,t)$ for $(x,t)\in{\bf R}^N\times(0,T)$.  
Let $\eta\in C^\infty([0,1])$ be such that 
$$
0\leq \eta\leq 1\quad\mbox{in}\quad[0,1],
\qquad 
\eta=1\quad\mbox{in}\quad[0,1/2],
\qquad
\eta=0\quad\mbox{in}\quad[3/4,1].
$$
Set $\eta_T(t):=\eta(t/T)$.
Then, for any integer $s > p/(p-1)$ and any nonnegative function $\zeta\in C^\infty_0({\bf R}^N)$, 
it follows from Definition~{\ref{Definition:1.1}} and the Young inequality that
\begin{equation}
\label{eq:2.2}
\begin{split}
 & \int_0^T\int_{{\bf R}^N}  u_z^p\zeta(x)^s\eta_T(t)^s\,dx\,dt + \int_0^T \int_{{\bf R}^N} \zeta(x)^s\eta_T(t)^s\,d\mu_z(x)\,dt\\
 & =\int_0^T\int_{{\bf R}^N} u_z\big(-\partial_t[\zeta(x)^s\eta_T(t)^s] + (-\Delta)^\frac{\theta}{2}[\zeta(x)^s\eta_T(t)^s]\big)\,dx\,dt\\
 & \le s\int_0^T\int_{{\bf R}^N} u_z\zeta(x)^s\eta_T(t)^{s-1}|\partial_t\eta_T|\,dx\,dt\\
 & \qquad\qquad
 +s\int_0^T\int_{{\bf R}^N} u_z\eta_T(t)^s\zeta(x)^{s-1}(-\Delta)^\frac{\theta}{2}\zeta(x)\,dx\,dt\\
 & \le \frac{1}{2}\int_0^T\int_{{\bf R}^N} u_z^p\zeta(x)^s\eta_T(t)^s\,dx\,dt
 +C\int_0^T\int_{{\bf R}^N}\zeta(x)^s\eta_T(t)^{s-\frac{p}{p-1}}|\partial_t\eta_T|^{\frac{p}{p-1}}\,dx\,dt\\
 & \qquad\qquad\qquad\qquad\qquad\qquad
 +C\int_0^T\int_{{\bf R}^N}\eta_T(t)^s\zeta^{s-\frac{p}{p-1}}|(-\Delta)^\frac{\theta}{2}\zeta(x)|^{\frac{p}{p-1}}\,dx\,dt.
\end{split}
\end{equation}
Here we also used the inequality 
$(-\Delta)^{\theta/2}\zeta^s \le s\zeta^{s-1} (-\Delta)^{\theta/2}\zeta$ 
(see \cite[Appendix]{FK12} and \cite[Proposition 3.3]{Ju05}).  
Since $s  > p/(p-1)$, $\eta_T=1$ on $[0,T/2]$ and $|\partial_t\eta_T|\le CT^{-1}$ on $[0,T]$, 
by \eqref{eq:2.2} we have
\begin{equation*}
\begin{split}
 & \int_0^{T/2} \int_{{\bf R}^N} \zeta(x)^s\,d\mu_z(x)\,dt\\
 & \le CT^{-\frac{p}{p-1}}\int_0^T\int_{{\bf R}^N}\zeta(x)^s\,dx\,dt
 +C\int_0^T\int_{{\bf R}^N}\zeta(x)^{s-\frac{p}{p-1}}|(-\Delta)^\frac{\theta}{2}\zeta(x)|^{\frac{p}{p-1}}\,dx\,dt.
\end{split}
\end{equation*}
This implies \eqref{eq:2.1}, and the proof is complete.
$\Box$
\vspace{5pt}
\newline
{\bf Proof of (\ref{eq:1.4}).}
Let $u$ be a  solution to \eqref{eq:1.1} in ${\bf R}^N\times[0,T)$ for some $T>0$. 
Let $z\in {\bf R}^N$ and $0<\sigma\le T^{1/\theta}$. 
Let $\zeta\in C_0^\infty({\bf R}^N)$ be such that 
$$
0\le \zeta\le 1\quad\mbox{in}\quad{\bf R}^N,
\qquad
\zeta=1\quad\mbox{in}\quad B(0,1/2),
\qquad
\mbox{supp}\,\zeta\subset B(0,1).
$$
Set $\zeta_\sigma(x):=\zeta(\sigma^{-1}(x-z))$ for $x\in{\bf R}^N$. 
By Lemma~\ref{Lemma:2.1} we have
\begin{equation}
\label{eq:2.3}
\mu(B(z,2^{-1}\sigma))
\leq C \sigma^{N-\frac{\theta p}{p-1}}
+ C\int_{B(z,\sigma)} 
|(-\Delta)^\frac{\theta}{2} \zeta_\sigma|^\frac{p}{p-1}
 \, dx.
\end{equation}
On the other hand, it follows that
\begin{equation}
\label{eq:2.4}
(-\Delta)^{\frac{\theta}{2}}\zeta_\sigma(x) 
= \sigma^{-\theta}[(-\Delta)^{\frac{\theta}{2}}\zeta](\sigma^{-1}(x-z)),
\end{equation}
which implies that 
\begin{equation}
\label{eq:2.5}
\begin{split}
\int_{B(z,\sigma)} 
|(-\Delta)^\frac{\theta}{2} \zeta_\sigma|^\frac{p}{p-1}\, dx
 & =\sigma^{-\frac{\theta p}{p-1}}\int_{B(0,\sigma)}
|[(-\Delta)^{\frac{\theta}{2}}\zeta](\sigma^{-1}x)|^{\frac{p}{p-1}}\,dx\\
 & =\sigma^{N-\frac{\theta p}{p-1}}\int_{B(0,1)}
|[(-\Delta)^{\frac{\theta}{2}}\zeta](x)|^{\frac{p}{p-1}}\,dx\le C\sigma^{N-\frac{\theta p}{p-1}}.
\end{split}
\end{equation}
By \eqref{eq:2.3} and \eqref{eq:2.5} we obtain 
$$
\sup_{z\in{\bf R}^N}\mu(B(z,2^{-1}\sigma))\le C\sigma^{N-\frac{\theta p}{p-1}}
$$
for $0<\sigma\le T^{1/\theta}$. 
By \cite[Lemma~2.1]{IS16} we find $M\in{\bf N}$ depending only on $N$ such that
\begin{equation}
\label{eq:2.6}
\sup_{z\in{\bf R}^N} \mu(B(z,\sigma))\le M \sup_{z\in{\bf R}^N}\mu(B(z,2^{-1}\sigma))\le CM\sigma^{N-\frac{\theta p}{p-1}}
\end{equation}
for $0<\sigma\le T^{1/\theta}$. 
Thus \eqref{eq:1.4} follows.
$\Box$
\vspace{5pt}

In the rest of this section we prove \eqref{eq:1.5} 
to complete the proof of Theorem~\ref{Theorem:1.1}. 
For this aim, 
it suffices to consider the case of 
\begin{equation}
\label{eq:2.7}
0<\theta<\min\{N,2\}\qquad\mbox{and}\qquad p=p_*\equiv N/(N-\theta). 
\end{equation}
Let $T>0$ and $0<\rho \le 32^{-1}T^{1/\theta}$. 
For $r>0$, set
\begin{equation}
\label{eq:2.8}
D_{r,\rho}:=B(0,r^{-1}T^\frac{1}{\theta})\setminus \overline{B(0, r\rho)},
\qquad
c_\rho:=\biggl(\displaystyle{\log\frac{T^\frac{1}{\theta}}{\rho}}\biggr)^{-1}.
\end{equation}
Let $f\in C^{\infty}_0(B(0,T^{1/\theta}))$ satisfy $f\ge0$ and 
\begin{equation*}
	f(x)\left\{
	\begin{aligned}
	&= 		c_\rho |x|^{-\theta},\qquad
	&& x\in D_{4,\rho}, \\
	&\le	c_\rho |x|^{-\theta}, 
	&& x\in D_{2,\rho}\setminus\overline{D_{4,\rho}},\\
	&=		0,
	&& \rm{otherwise}.
	\end{aligned}\right.
\end{equation*}
Consider the following fractional Poisson equation 
\begin{equation}
\label{eq:2.9}
	\left\{
	\begin{aligned}
	&(-\Delta)^\frac{\theta}{2} \psi(x)= f(x),\qquad
	&& x\in B(0,T^\frac{1}{\theta}), \\
	&\psi(x)= 0, 
	&& x\in {\bf R}^N\setminus B(0,T^\frac{1}{\theta}).
	\end{aligned}\right.
\end{equation}
It follows from \cite{Bu16} that 
the solution $\psi$ is represented as 
\begin{equation}\label{eq:2.10}
\psi(x) = \left\{
	\begin{aligned}
	&\int_{D_{2,\rho}} f(y)\Gamma(x,y)\,dy,
	&&x\in B(0,T^\frac{1}{\theta}), \\
	&0, && x\in {\bf R}^N\setminus B(0,T^\frac{1}{\theta}).
	\end{aligned}
\right.
\end{equation}
Here 
\begin{equation}
\label{eq:2.11}
\Gamma(x,y) := \kappa|x-y|^{\theta-N} \int_0^{r_0(x,y)} \frac{\tau^{\frac{\theta}{2}-1}}{(\tau+1)^\frac{N}{2}}\,d\tau,
\end{equation}
where $\kappa=\kappa(N,\theta)>0$ and
\begin{equation}
\label{eq:2.12}
r_0(x,y):=\frac{(T^\frac{2}{\theta}-|x|^2)(T^\frac{2}{\theta}-|y|^2)}{T^\frac{2}{\theta}|x-y|^2}.
\end{equation}
Then, by \cite[Proposition~1.1]{RS14} we see that  
$\psi\in C^{\theta/2}({\bf R}^N)$. 
Furthermore, we have 
\begin{lemma}
\label{Lemma:2.2}
Assume the same conditions as in Theorem~{\rm\ref{Theorem:1.1}}. 
Let $0<\theta<\min\{N,2\}$, $p=p_*$ and $\psi$ be as in the above. 
Then there exists $C=C(N,\theta)>0$ such that
\begin{eqnarray}
\label{eq:2.13}
 & & \sup_{B(0,T^{1/\theta})}\psi \le C,\\
\label{eq:2.14}
 & & \inf_{B(0,\rho)} \psi \ge C^{-1},\\
\label{eq:2.15}
 & & \int_{B(0,T^{1/\theta})} |(-\Delta)^\frac{\theta}{2} \psi|^\frac{p}{p-1}\, dx 
 \le C c_\rho^{\frac{1}{p-1}},\\
\label{eq:2.16}
 & & \int_{{\bf R}^N} \psi^\frac{p}{p-1} \, dx=\int_{B(0,T^{1/\theta})} \psi^\frac{p}{p-1} \, dx 
 \le C  c_\rho^{\frac{1}{p-1}}T^\frac{p}{p-1},
\end{eqnarray}
for $0<\rho \le 32^{-1}T^{1/\theta}$.
\end{lemma}
{\bf Proof.}
We prove \eqref{eq:2.13}. 
Let $x\in B(0,T^{1/\theta})$.
It follows from \eqref{eq:2.7} that 
\begin{equation}
\label{eq:2.17}
\int_0^{\infty} \frac{\tau^{\frac{\theta}{2}-1}}{(\tau+1)^\frac{N}{2}}\,d\tau <\infty,
\qquad
\Gamma(x,y)\le C|x-y|^{\theta-N}.
\end{equation}
Then, by \eqref{eq:2.8} we have
\begin{equation*}
\begin{split}
 & \int_{D_{2,\rho}\cap B(x,2\rho)} f(y) \Gamma(x,y) \, dy
\le Cc_\rho\int_{D_{2,\rho}\cap B(x,2\rho)} |y|^{-\theta} |x-y|^{\theta-N} \, dy\\
 & \qquad\quad
\le Cc_\rho\rho^{-\theta} \int_{B(x,2\rho)} |x-y|^{\theta-N}\,dy
\le Cc_\rho
\le C.
\end{split}
\end{equation*}
Similarly, by the H\"older inequality and \eqref{eq:2.17}
we have
\begin{equation*}
\begin{split}
 & \int_{D_{2,\rho}\setminus B(x,2\rho)} f(y) \Gamma(x,y) \, dy\le Cc_\rho\int_{D_{2,\rho}\setminus B(x,2\rho)}|y|^{-\theta}|x-y|^{\theta-N}\\
 &\le 	Cc_\rho\biggl( \int_{D_{2,\rho}} |y|^{-N} \,dy \biggr)^{\frac{\theta}{N}} 
		\biggl( \int_{B(x,2T^{1/\theta})\setminus B(x,2\rho)} |x-y|^{-N} \,dy \biggr)^{\frac{N-\theta}{N}}
		\le Cc_\rho\log\frac{T^\frac{1}{\theta}}{\rho}=C.
\end{split}
\end{equation*}
These imply that 
$$
\psi(x)=\int_{D_{2,\rho}\cap B(x,2\rho)} f(y) \Gamma(x,y) \, dy+
\int_{D_{2,\rho}\setminus B(x,2\rho)} f(y) \Gamma(x,y) \, dy\le C 
$$
for $x\in B(0,T^{1/\theta})$. Thus \eqref{eq:2.13} follows. 

We prove \eqref{eq:2.14}. 
Since $\rho \le 32^{-1}T^{1/\theta}$, 
by \eqref{eq:2.12} we find $c>0$ such that  
$$
r_0(x,y)
\ge \frac{(T^\frac{2}{\theta}-\rho^2)
(T^\frac{2}{\theta} -16^{-1}T^\frac{2}{\theta})}{T^\frac{2}{\theta}(4^{-1}T^\frac{1}{\theta}+\rho)^2} 
\ge c>0,
\qquad
x\in B(0,\rho),\,\,\, y\in D_{4,\rho}.
$$
It follows that 
$|x-y|\le |x|+|y| \le (5/4)|y|$ for $x\in B(0,\rho)$ and $y\in D_{4,\rho}$.
 Then, by \eqref{eq:2.10} and \eqref{eq:2.11} we have
\begin{equation*}
\begin{split}
	\psi(x) 
	&\ge \kappa c_\rho
	\int_0^c \frac{\tau^{\frac{\theta}{2}-1}}{(\tau+1)^\frac{N}{2}}	\,d\tau
	\int_{D_{4,\rho}} |y|^{-\theta}|x-y|^{\theta-N} \,dy\\
	&\geq C^{-1}c_\rho
	\int_{D_{4,\rho}} |y|^{-\theta} \left( \frac{5}{4}|y|\right)^{\theta-N} \,dy
	\geq C^{-1}c_\rho\log\frac{T^\frac{1}{\theta}}{\rho}
	=C^{-1}
\end{split}
\end{equation*}
for $x\in B(0,\rho)$. This implies \eqref{eq:2.14}.

We prove \eqref{eq:2.15} and \eqref{eq:2.16}. 
Since $p=p_*=N/(N-\theta)$, by \eqref{eq:2.9} we have
\[
\begin{split}
	 & \int_{B(0,T^{1/\theta})} |(-\Delta)^\frac{\theta}{2} \psi|^\frac{p}{p-1}\, dx
	=\int_{B(0,T^{1/\theta})} f(x)^\frac{p}{p-1}\, dx\\
	 & \le c_\rho^{\frac{p}{p-1}}\int_{D_{2,\rho}}|x|^{-\frac{\theta p}{p-1}}\,dx
	=c_\rho^{\frac{p}{p-1}}\int_{D_{2,\rho}}|x|^{-N}\,dx
	\leq Cc_\rho^{\frac{p}{p-1}}\log\frac{T^\frac{1}{\theta}}{\rho}
	=C c_\rho^{\frac{1}{p-1}},
\end{split}
\]
which implies \eqref{eq:2.15}.
Furthermore, 
by \eqref{eq:2.10} and \eqref{eq:2.17} 
we apply the Young inequality to obtain 
\[
\begin{split}
\int_{B(0,T^{1/\theta})} \psi^\frac{p}{p-1} \, dx
& \leq C c_\rho^{\frac{p}{p-1}}\int_{B(0,T^{1/\theta})} \biggl( \int_{D_{2,\rho}} |y|^{-\theta} |x-y|^{\theta-N} \, dy \biggr)^\frac{p}{p-1}\, dx\\
& \le C c_\rho^{\frac{p}{p-1}}\int_{{\bf R}^N} \biggl( \int_{{\bf R}^N} |y|^{-\theta} 
 \chi_{D_{2,\rho}}(y) |x-y|^{\theta-N}\chi_{B(0,2T^{1/\theta})}(x-y) \, dy \biggr)^{\frac{N}{\theta}}\, dx\\
&\leq	Cc_\rho^{\frac{p}{p-1}}
 \biggl(\int_{{\bf R}^N} |x|^{-N} \chi_{D_{2,\rho}}(x) \, dx \biggr)
 \biggl(\int_{{\bf R}^N} |x|^{\theta-N} \chi_{B(0,2T^{1/\theta})}(x) \, dx \biggr)^\frac{N}{\theta}\\
&\leq Cc_\rho^{\frac{p}{p-1}}T^\frac{N}{\theta}\log\frac{T^\frac{1}{\theta}}{\rho}
=Cc_\rho^{\frac{1}{p-1}}T^{\frac{p}{p-1}}.
\end{split}
\]
Thus \eqref{eq:2.16} holds, and the proof of Lemma~\ref{Lemma:2.2} is complete. 
$\Box$
\begin{lemma}
\label{Lemma:2.3}
Let $u$ be a solution to problem~\eqref{eq:1.1} in ${\bf R}^N\times[0,T)$ for some $T>0$. 
Let $p=p_*$ and $\psi$ be as in the above. 
Let $s$ be a sufficiently large integer such that $s > p/(p-1)$. 
Then \eqref{eq:2.1} holds with $\zeta$ replaced by $\psi$. 
\end{lemma}
{\bf Proof.}
Let $z\in{\bf R}^N$. 
Let $\eta\in C^\infty_0({\bf R}^N)$ be such that
$$
\eta\ge 0\quad\mbox{in}\quad{\bf R}^N,\qquad
\mbox{supp}\,\eta\subset B(0,1),\qquad
\int_{{\bf R}^N} \eta(x) \, dx =1.
$$
For any $\epsilon>0$, set 
\begin{equation}
\label{eq:2.18}
\eta_\epsilon(x) := (\epsilon T^\frac{1}{\theta})^{-N} \eta\biggl(\frac{x}{\epsilon T^\frac{1}{\theta}}\biggr),
\qquad x\in{\bf R}^N. 
\end{equation}
By \eqref{eq:2.10} we have 
$$
\psi*\eta_\epsilon \in C^\infty_0({\bf R}^N),\qquad
\mbox{supp}\,(\psi*\eta_\epsilon)\subset 
B(0,(1+\epsilon)T^{1/\theta}).
$$ 
Let $s$ be an integer such that $s > p/(p-1)=N/\theta$. 
Then it follows from Lemma~\ref{Lemma:2.1} that 
\begin{equation}
\label{eq:2.19}
\begin{split}
 & \int_{B(0,(1+\epsilon)T^{1/\theta})} (\psi*\eta_\epsilon)^s \, d\mu_z(x)
\leq C T^{-\frac{p}{p-1}} \int_{B(0,(1+\epsilon)T^{1/\theta})} (\psi*\eta_\epsilon)^s  \, dx \\
 & \qquad\qquad
+ C\int_{B(0,(1+\epsilon)T^{1/\theta})} 
|(-\Delta)^\frac{\theta}{2} (\psi*\eta_\epsilon)|^\frac{p}{p-1}(\psi*\eta_\epsilon)^{s-\frac{p}{p-1}} \, dx.
\end{split}
\end{equation}
Furthermore, it is easy to see that 
\begin{equation}
\label{eq:2.20}
\begin{split}
&\lim_{\epsilon\to0} \int_{B(0,(1+\epsilon)T^{1/\theta})} (\psi*\eta_\epsilon)^s \, d\mu_z(x)
= \int_{B(0,T^{1/\theta})} \psi^s \, d\mu_z(x),  \\
&\lim_{\epsilon\to0} \int_{B(0,(1+\epsilon)T^{1/\theta})} (\psi*\eta_\epsilon)^s  \, dx 
= \int_{B(0,T^{1/\theta})} \psi^s  \, dx.
\end{split} 
\end{equation}

We prove 
\begin{equation}
\label{eq:2.21}
\begin{split}
 & \lim_{\epsilon\to0}\int_{B(0,(1+\epsilon)T^{1/\theta})} 
|(-\Delta)^\frac{\theta}{2} (\psi*\eta_\epsilon)|^\frac{p}{p-1}(\psi*\eta_\epsilon)^{s-\frac{p}{p-1}} \, dx\\
 & =\int_{B(0,T^{1/\theta})} |(-\Delta)^\frac{\theta}{2}\psi|^\frac{p}{p-1}\psi^{s-\frac{p}{p-1}} \, dx.
\end{split}
\end{equation}
It follows that 
\begin{equation}
\label{eq:2.22}
\begin{split}
(-\Delta)^\frac{\theta}{2} (\psi*\eta_\epsilon)(x) 
&=		\mathcal{F}^{-1}[|\xi|^\theta \widehat{\psi*\eta_\epsilon}](x)
 =\mathcal{F}^{-1}[|\xi|^\theta \hat{\psi}\hat{\eta_\epsilon}](x)\\ 
&=		\mathcal{F}^{-1}[|\xi|^\theta\hat{\psi}] * \mathcal{F}^{-1}[\hat{\eta_\epsilon}](x)
=([(-\Delta)^\frac{\theta}{2}\psi]*\eta_\epsilon)(x)
\end{split}
\end{equation}
for $x\in{\bf R}^N$. 
For any $x\in B(0,(1-\epsilon)T^{1/\theta})$ and $y\in B(0,\epsilon T^{1/\theta})$, 
since $x-y\in B(0,T^{1/\theta})$, 
by \eqref{eq:2.9} and \eqref{eq:2.22} we have
$$
(-\Delta)^\frac{\theta}{2} (\psi*\eta_\epsilon)(x)
=\int_{B(0,\epsilon T^{1/\theta})}((-\Delta)^\frac{\theta}{2}\psi)(x-y)\eta_\epsilon(y)\,dy
=(f*\eta_\epsilon)(x).
$$
Then 
\begin{equation}
\label{eq:2.23}
\begin{split}
I_\epsilon := & 
\int_{B(0,(1-\epsilon)T^{1/\theta})}|(-\Delta)^\frac{\theta}{2} 
 (\psi*\eta_\epsilon)|^\frac{p}{p-1}(\psi*\eta_\epsilon)^{s-\frac{p}{p-1}} \, dx\\
= & \,\int_{B(0,(1-\epsilon)T^{1/\theta})} (f*\eta_\epsilon)^\frac{p}{p-1}(\psi*\eta_\epsilon)^{s-\frac{p}{p-1}} \, dx
\to\int_{B(0,T^{1/\theta})} f^\frac{p}{p-1}\psi^{s-\frac{p}{p-1}} \, dx
\end{split}
\end{equation}
as $\epsilon\to 0$. 

Similarly to \eqref{eq:2.4} and \eqref{eq:2.22}, 
it follows from \eqref{eq:2.18} that 
\begin{equation*}
\begin{split}
(-\Delta)^\frac{\theta}{2}\eta_\epsilon(x)
 & =(\epsilon T^\frac{1}{\theta})^{-N-\theta} [(-\Delta)^\frac{\theta}{2}\eta]\biggl(\frac{x}{\epsilon T^\frac{1}{\theta}}\biggr),\\
(-\Delta)^\frac{\theta}{2}(\psi*\eta_\epsilon)(x) & = (\psi*[(-\Delta)^\frac{\theta}{2}\eta_\epsilon])(x)\\
 & =(\epsilon T^\frac{1}{\theta})^{-N-\theta}
\int_{{\bf R}^N}\psi(x-y)[(-\Delta)^\frac{\theta}{2}\eta]\left(\frac{y}{\epsilon T^\frac{1}{\theta}}\right)\,dy\\
 & =(\epsilon T^\frac{1}{\theta})^{-\theta}
\int_{{\bf R}^N}\psi(x-\epsilon T^{\frac{1}{\theta}}y)(-\Delta)^\frac{\theta}{2}\eta(y)\,dy.
\end{split}
\end{equation*}
Then 
\begin{equation}
\label{eq:2.24}
\begin{split}
J_\epsilon
 := & \int_{E_\epsilon}|(-\Delta)^\frac{\theta}{2} 
 (\psi*\eta_\epsilon)|^\frac{p}{p-1}(\psi*\eta_\epsilon)^{s-\frac{p}{p-1}}\,dx\\
 = &\,(\epsilon T^\frac{1}{\theta})^{-\frac{\theta p}{p-1}} \int_{E_\epsilon}  \biggl|\int_{{\bf R}^N} \psi(x-\epsilon T^\frac{1}{\theta}y)(-\Delta)^\frac{\theta}{2}\eta (y) \, dy\biggr|^\frac{p}{p-1} (\psi*\eta_\epsilon)(x)^{s-\frac{p}{p-1}}\,dx,\\
\end{split}
\end{equation}
where $E_\epsilon:=B(0, (1+\epsilon)T^{1/\theta})\setminus B(0, (1-\epsilon)T^{1/\theta})$. 
Since $\eta\in C_0^\infty({\bf R}^N)\subset H^{\theta,p}({\bf R}^N)$ 
(see e.g. \cite[Theorem~7.38]{A}), 
by the H\"{o}lder inequality and \eqref{eq:2.16} we see that
\begin{equation}
\label{eq:2.25}
\begin{split}
&		\biggl|\int_{{\bf R}^N} \psi(x-\epsilon T^\frac{1}{\theta}y)(-\Delta)^\frac{\theta}{2}\eta (y) \, dy\biggr|^\frac{p}{p-1}\\
&\le \biggl(\int_{{\bf R}^N} \psi(x-\epsilon T^\frac{1}{\theta}y)^\frac{p}{p-1} \, dy\biggr)\biggl( \int_{{\bf R}^N} |(-\Delta)^\frac{\theta}{2}\eta(y)|^p \, dy\biggr)^\frac{1}{p-1}\\
&\le 	C(\epsilon T^\frac{1}{\theta})^{-N} 
\int_{{\bf R}^N} \psi^\frac{p}{p-1} \, dy
\le C(\epsilon T^\frac{1}{\theta})^{-N}
\end{split}
\end{equation}
for $x\in E_\epsilon$. 
On the other hand, it follows that
\begin{equation}
\label{eq:2.26}
\begin{split}
\int_{E_\epsilon} (\psi*\eta_\epsilon)(x)^{s-\frac{p}{p-1}}\,dx
&=		\int_{E_\epsilon} \biggl(\int_{{\bf R}^N} \psi(x-y) (\epsilon T^\frac{1}{\theta})^{-N}\eta\biggl(\frac{y}{\epsilon T^\frac{1}{\theta}}\biggr) \, dy \biggr)^{s-\frac{p}{p-1}}\,dx\\
&\le		\int_{E_\epsilon} \biggl(\int_{B(0,1)} \psi(x-\epsilon T^\frac{1}{\theta}y) \eta(y) \, dy \biggr)^{s-\frac{p}{p-1}}\,dx
\end{split}
\end{equation}
for $x\in E_\epsilon$. 
Since 
$$
\mbox{dist}(x-\epsilon T^{\frac{1}{\theta}}y, \partial B(0,T^{\frac{1}{\theta}}) )
<2\epsilon T^{\frac{1}{\theta}}
$$ 
for $x\in E_\epsilon$ and $y\in B(0,1)$, 
recalling that $\psi = 0$ in ${\bf R}^N\setminus B(0,T^{1/\theta})$ and $\psi\in C^{\theta/2}({\bf R}^N)$, 
we observe that
$$ 
0\le \psi(x-\epsilon T^\frac{1}{\theta}y)\le (2\epsilon T^\frac{1}{\theta})^\frac{\theta}{2} \|\psi\|_{C^{\theta/2}({\bf R}^N)}
$$
for $x\in E_\epsilon$ and $y\in B(0,1)$.
This together with \eqref{eq:2.26} implies that
\begin{equation}
\label{eq:2.27}
\int_{E_\epsilon} (\psi*\eta_\epsilon)(x)^{s-\frac{p}{p-1}}\,dx
\le C \epsilon^{\frac{\theta}{2}\left(s-\frac{p}{p-1}\right)+1} T^{\frac{1}{2}\left(s-\frac{p}{p-1}\right)+\frac{N}{\theta}}
\|\psi\|_{C^{\theta/2}({\bf R}^N)}.
\end{equation}
Therefore, taking a sufficiently large integer $s > p/(p-1)$ if necessary, 
we deduce from \eqref{eq:2.24}, \eqref{eq:2.25} and \eqref{eq:2.27} that
\begin{equation}
\label{eq:2.28}
\begin{split}
\limsup_{\epsilon\to 0}J_\epsilon
& \le	C
T^{-\frac{p}{p-1}-\frac{N}{\theta}+\frac{1}{2}\left(s-\frac{p}{p-1}\right)+\frac{N}{\theta}}
\|\psi\|_{C^{\theta/2}({\bf R}^N)}\lim_{\epsilon\to 0}
\epsilon^{-\frac{\theta p}{p-1}-N+\frac{\theta}{2}\left(s-\frac{p}{p-1}\right)+1}=0.
\end{split}
\end{equation}
Then, by \eqref{eq:2.23} and \eqref{eq:2.28} 
we have
\begin{equation*}
\begin{split}
 & \lim_{\epsilon\to0}\int_{B(0,(1+\epsilon)T^{1/\theta})} 
|(-\Delta)^\frac{\theta}{2} (\psi*\eta_\epsilon)|^\frac{p}{p-1}(\psi*\eta_\epsilon)^{s-\frac{p}{p-1}} \, dx\\
 & =\lim_{\epsilon\to0}I_\epsilon+\lim_{\epsilon\to0}J_\epsilon
=\int_{B(0,T^{1/\theta})} f^\frac{p}{p-1}\psi^{s-\frac{p}{p-1}} \, dx,
\end{split}
\end{equation*}
which implies \eqref{eq:2.21}. 
Combining \eqref{eq:2.19}, \eqref{eq:2.20} and \eqref{eq:2.21}, 
we obtain \eqref{eq:2.1} with $\zeta$ replaced by $\psi$. 
Thus Lemma~\ref{Lemma:2.3} follows.
$\Box$
\vspace{5pt}

Now we are ready to complete the proof of Theorem~\ref{Theorem:1.1}. 
\vspace{5pt}
\newline
{\bf Proof of Theorem~\ref{Theorem:1.1}.}
It remains to prove \eqref{eq:1.5}. 
Let $u$ be a solution to \eqref{eq:1.1} in ${\bf R}^N \times[0,T)$ and $p=p_*$. 
Let $z\in{\bf R}^N$, $0<\rho \le 32^{-1} T^{1/\theta}$ and $s$ be as in Lemma~\ref{Lemma:2.3}. 
It follows from Lemmas~\ref{Lemma:2.2} and \ref{Lemma:2.3} that
\begin{equation}
\label{eq:2.29}
\begin{split}
\int_{{\bf R}^N} \psi^s \, d\mu_z(x) 
 & \leq C T^{-\frac{p}{p-1}} \int_{{\bf R}^N} \psi^s  \, dx 
+ C\int_{{\bf R}^N} |(-\Delta)^\frac{\theta}{2} \psi|^\frac{p}{p-1}\psi^{s-\frac{p}{p-1}} \, dx\\
 &  \leq C T^{-\frac{p}{p-1}} \int_{{\bf R}^N} \psi^{\frac{p}{p-1}}  \, dx 
+ C \int_{B(0,T^{\frac{1}{\theta}})} |(-\Delta)^\frac{\theta}{2} \psi|^\frac{p}{p-1} \, dx\\
 & \le Cc_\rho^{\frac{1}{p-1}}=C\biggl(\log\frac{T^\frac{1}{\theta}}{\rho}\biggr)^{-\frac{1}{p-1}}.
\end{split}
\end{equation}
On the other hand, by \eqref{eq:2.10} and \eqref{eq:2.14} we have 
\begin{equation}
\label{eq:2.30}
\begin{split}
\int_{{\bf R}^N} \psi^{s} \, d\mu_z(x) 
\ge 
\int_{B(0,\rho)}
\psi^{s} \, d\mu_z(x)
\ge C\mu(B(z,\rho)),\quad z\in{\bf R}^N.
\end{split}
\end{equation}
By \eqref{eq:2.29} and \eqref{eq:2.30} we obtain  
$$
\mu(B(z,\rho))
\le C\biggl(\log\frac{T^\frac{1}{\theta}}{\rho}\biggr)^{-\frac{N}{\theta}+1} 
\le C\left[\log\biggl( e +\frac{T^\frac{1}{\theta}}{\rho}\biggr)\right]^{-\frac{N}{\theta}+1}
$$
for $0<\rho\le32^{-1}T^{1/\theta}$.
Then, similarly to \eqref{eq:2.6}, we deduce that
$$
\sup_{z\in{\bf R}^N}\,\mu(B(z,\sigma))
\le C\left[\log\biggl( e +\frac{T^\frac{1}{\theta}}{\sigma}\biggr)\right]^{-\frac{N}{\theta}+1}
$$
for $0<\sigma\le T^{1/\theta}$.
Thus \eqref{eq:1.5} holds, and the proof of Theorem~\ref{Theorem:1.1} is complete. 
$\Box$
\vspace{3pt}

We prove Corollaries~\ref{Corollary:1.1} and \ref{Corollary:1.2}. 
\vspace{3pt}
\newline
{\bf Proof of Corollary~\ref{Corollary:1.1}.}
Under the assumptions of Corollary~\ref{Corollary:1.1}, 
it follows that 
\begin{equation*}
\begin{split}
 & \sup_{x\in{\bf R}^N}\int_{B(x,\sigma)}\mu(y)\,dy\ge C\gamma\sigma^{N-\frac{\theta p}{p-1}}
  \qquad\,\,\,\,\mbox{if}\quad p>p_*,\\
 & \sup_{x\in{\bf R}^N}\int_{B(x,\sigma)}\mu(y)\,dy\ge C\gamma|\log\sigma|^{-\frac{N}{\theta}+1}
 \quad\mbox{if}\quad p=p_*,
\end{split}
\end{equation*}
for sufficiently small $\sigma>0$. 
These together with Theorem~\ref{Theorem:1.1} imply that 
problem~\eqref{eq:1.1} possesses no local-in-time solutions if $\gamma$ is sufficiently large. 
Thus Corollary~\ref{Corollary:1.1} follows.
$\Box$
\vspace{3pt}
\newline
{\bf Proof of Corollary~\ref{Corollary:1.2}.}
Suppose that problem~\eqref{eq:1.1} possesses a global-in-time solution. 
Let $1<p\leq p_*$. 
Then it follows from Theorem~\ref{Theorem:1.1} that  
\begin{equation*}
\begin{split}
 & \sup_{x\in{\bf R}^N}\mu(B(x,T^\frac{1}{\theta}))\le \gamma\,T^{\frac{N}{\theta}-\frac{p}{p-1}}
 \qquad\qquad\qquad\quad\,\,\,\mbox{if}\quad 1<p<p_*,\\
 & \sup_{x\in{\bf R}^N}\mu(B(x,T^\frac{1}{2\theta}))\le \gamma\,
 \left[\log\left(e+T^{\frac{1}{2\theta}}\right)\right]^{-\frac{N}{\theta}+1}
 \quad\mbox{if}\quad p=p_*,
\end{split}
\end{equation*}
for sufficiently large $T>0$. Letting $T\to\infty$, we see that  
$\mu({\bf R}^N)=0$, which contradicts $\mu\not\equiv 0$ in ${\bf R}^N$. 
Thus assertion~(1) follows.

Assume that there exists $R>0$ such that 
$\mu(x)\ge\gamma|x|^{-\theta p/(p-1)}$ 
for almost all $x\in{\bf R}^N\setminus B(0,R)$. 
Then we have
\begin{equation}
\label{eq:2.31}
\sup_{x\in{\bf R}^N}\mu(B(x,\sigma))\ge C\gamma \sigma^{N-\frac{\theta p}{p-1}}
\end{equation}
for $\sigma\ge 2R$. If $\gamma$ is sufficiently large, then \eqref{eq:2.31} 
contradicts \eqref{eq:1.4} for sufficiently large $\sigma$. 
This means that problem~\eqref{eq:1.1} possesses no global-in-time solutions. 
Thus assertion~(2) follows, and the proof of Corollary~\ref{Corollary:1.2} is complete.
$\Box$
\section{Proofs of Theorems~\ref{Theorem:1.2} and \ref{Theorem:1.3}}
Consider the following integral equation
\begin{equation}
\tag{I}
u(x,t)
=\int_0^t \int_{{\bf R}^N}G(x-y,t-s)\,d\mu(y)\,ds 
+\int_0^t\int_{{\bf R}^N}G(x-y,t-s)u(y,s)^p\,dy\,ds,
\end{equation}
where $N\ge 1$, $0<\theta\le 2$ and $p>1$. 
Here $G=G(x,t)$ is the fundamental solution to 
\begin{equation}
\label{eq:3.1}
\partial_t v +(-\Delta)^{\frac{\theta}{2}} v = 0 \qquad \mbox{in} \qquad {\bf R}^N\times(0,\infty).
\end{equation}
In particular, in the case of $\theta=2$, 
\begin{equation}
\label{eq:3.2}
G(x,t)=(4\pi t)^{-\frac{N}{2}}\exp\left(-\frac{|x|^2}{4t}\right). 
\end{equation}
We obtain sufficient conditions on the solvability of integral equation~(I)
and prove Theorems~\ref{Theorem:1.2} and \ref{Theorem:1.3}. 
\subsection{Preliminaries}
We formulate the definition of solutions to integral equation~(I). 
\begin{definition}
\label{Definition:3.1}
Let $u$ be a nonnegative measurable function in ${\bf R}^N\times(0,T)$, 
where $0<T\le\infty$. 
We say that $u$ is a solution to integral equation~{\rm (I)} 
in ${\bf R}^N\times[0,T)$ if $u$ satisfies 
\begin{equation}
\label{eq:3.3}
\infty>u(x,t)
=\int_0^t \int_{{\bf R}^N}G(x-y,t-s)\,d\mu(y)\,ds 
+\int_0^t\int_{{\bf R}^N}G(x-y,t-s)u(y,s)^p\,dy\,ds
\end{equation}
for almost all $x\in{\bf R}^N$ and $0<t<T$. 
We say that $u$ is a supersolution to integral equation~{\rm (I)} 
if $u$ satisfies \eqref{eq:3.3} with ``\,$=$'' replaced by ``\,$\ge$'' 
for almost all $x\in{\bf R}^N$ and $0<t<T$.  
\end{definition}

We recall the following properties of the fundamental solution~$G$:  
\begin{align}
\notag
 & G\in C^\infty({\bf R}^N\times(0,\infty)),\quad
 G(x,t)=t^{-\frac{N}{\theta}}G(t^{-\frac{1}{\theta}}x,1),\\
\label{eq:3.4}
 & C^{-1}t^{-\frac{N}{\theta}}
 \left( 1+ \frac{|x|}{t^{1/\theta}} \right)^{-N-\theta} \leq 
 G(x,t)\le C t^{-\frac{N}{\theta}}
 \left( 1+ \frac{|x|}{t^{1/\theta}} \right)^{-N-\theta}\mbox{ if }\quad 0<\theta<2,\\
\label{eq:3.5}
 &  \int_{{\bf R}^N}G(x,t)\,dx=1, 
\end{align}
for $x$, $y\in{\bf R}^N$ and $t>s>0$ (see e.g. \cite{BSS03, HI18, Sugi75}).  
Furthermore, it follows that 
\begin{equation}
\label{eq:3.6}
\int_{{\bf R}^N} (-\Delta)^\frac{\theta}{2}\varphi(x)G(x,t) \, dx = \int_{{\bf R}^N} \varphi(x)(-\Delta)^\frac{\theta}{2}G(x,t) \, dx,
\qquad t>0,
\end{equation}
for $\varphi\in C^\infty_0({\bf R}^N)$. 
We show that a solution to integral equation~(I) in ${\bf R}^N\times[0,T)$ 
is a solution to problem~\eqref{eq:1.1} in ${\bf R}^N\times[0,T)$, where $0<T\le\infty$.
\begin{lemma}
\label{Lemma:3.1}
Let $u$ be a solution to integral equation~{\rm (I)} in ${\bf R}^N\times[0,T)$, 
where $0<T\le\infty$. Then 
$u\in L^p_{loc}({\bf R}^N\times[0,T))$ and 
\begin{equation}
\label{eq:3.7}
\int_0^T\int_{{\bf R}^N} u(-\partial_t \varphi + (-\Delta)^\frac{\theta}{2} \varphi )\,dx\,dt = 
\int_0^T\int_{{\bf R}^N}  u^p\varphi\,dx\,dt + \int_0^T \int_{{\bf R}^N}  \varphi \,d\mu(x)\,dt
\end{equation}
for $\varphi\in C^\infty_0({\bf R}^N\times[0,T))$.
\end{lemma}
{\bf Proof.}
It suffices to consider the case of $0<T<\infty$.
Let $u$ be a solution to integral equation~{\rm (I)} in ${\bf R}^N\times[0,T)$. 
Let $\epsilon\in(0,T/2)$ and $0<\theta<2$. 
By \eqref{eq:3.3} and \eqref{eq:3.4}
we find $t\in(T-\epsilon,T)$ such that 
\begin{equation*}
\begin{split}
\infty & >u(x,t)\ge \int_0^{T-2\epsilon}\int_{{\bf R}^N}G(x-y,t-s)\,[d\mu(y)+u(y,s)^p\,dy]\,ds\\
 & \ge C^{-1}T^{-\frac{N}{\theta}}\int_0^{T-2\epsilon}\int_{{\bf R}^N}
 \left(1+\frac{|x-y|}{\epsilon^{1/\theta}}\right)^{-N-\theta}
 \,[d\mu(y)+u(y,s)^p\,dy]\,ds
\end{split}
\end{equation*}
for almost all $x\in{\bf R}^N$. 
Since $\epsilon\in(0,T/2)$ is arbitrary, 
we see that $u\in L^p_{loc}({\bf R}^N\times[0,T))$. 
In the case of $\theta=2$, using \eqref{eq:3.2}, instead of \eqref{eq:3.4}, we 
deduce that $u\in L^p_{loc}({\bf R}^N\times[0,T))$.

Let $\varphi\in C^\infty_0({\bf R}^N\times[0,T))$.
It follows from \eqref{eq:3.1} and \eqref{eq:3.6} that
\begin{equation*}
\begin{split}
 &  \int_0^T\int_{{\bf R}^N}  \varphi(y,s)\,d\mu(y)\,ds\\
 & =\int_0^T\int_{{\bf R}^N} 
  \left( \int_s^T \int_{{\bf R}^N} (\partial_t+(-\Delta)^\frac{\theta}{2})\,
  G(x-y,t-s)\varphi(x,t)\,dx\,dt + \varphi(y,s) \right) \,d\mu(y)\,ds\\
  & =\int_0^T\int_{{\bf R}^N} \int_s^T \int_{{\bf R}^N} 
  G(x-y,t-s)(-\partial_t+(-\Delta)^\frac{\theta}{2})\,\varphi(x,t)\,dx\,dt\,d\mu(y)\,ds\\
 & =\int_0^T\int_{{\bf R}^N} \left(\int_0^t \int_{{\bf R}^N} 
 G(x-y,t-s)\,d\mu(y)\,ds\right)(-\partial_t+(-\Delta)^\frac{\theta}{2})\,\varphi(x,t)\,dx\,dt.
\end{split}
\end{equation*}
Similarly, we have
\[
\begin{aligned}
 &  \int_0^T\int_{{\bf R}^N} u^p \varphi \,dy\,ds\\
 & =\int_0^T\int_{{\bf R}^N} \left(\int_0^t \int_{{\bf R}^N} 
 G(x-y,t-s)u(y,s)^p\,dy\,ds\right)(-\partial_t+(-\Delta)^\frac{\theta}{2})\,\varphi(x,t)\,dx\,dt.
\end{aligned}
\]
Then 
\begin{equation*}
\begin{split}
 & \int_0^T\int_{{\bf R}^N}  u^p\varphi\,dx\,dt + \int_0^T \int_{{\bf R}^N}  \varphi \,d\mu(x)\,dt\\
 & =\int_0^T\int_{{\bf R}^N}\left(\int_0^t\int_{{\bf R}^N}G(x-y,t-s)
[u(y,s)^p\,dy+\,d\mu(y)]\,ds\right) (-\partial_t+(-\Delta)^\frac{\theta}{2})\varphi(x,t)\,dx\,dt\\
 & =\int_0^T\int_{{\bf R}^N}u(x,t)(-\partial_t+(-\Delta)^\frac{\theta}{2})\varphi(x,t)\,dx\,dt,
\end{split}
\end{equation*}
which implies \eqref{eq:3.7}. Thus Lemma~\ref{Lemma:3.1} follows. 
$\Box$

\subsection{Proofs of Theorems~\ref{Theorem:1.2} and \ref{Theorem:1.3}}
In this subsection we first apply the contraction mapping theorem in uniformly local Lebesgue spaces 
to integral equation~(I) and prove Theorem~\ref{Theorem:1.3}. 
Next we prove Theorem~\ref{Theorem:1.2}. 

Let $N\geq 1$, $0<\theta\le2$ and $1<p<p_*$. 
Assume the same conditions as in Theorem~\ref{Theorem:1.3} and set $T=\sigma^\theta$. 
Let $u_T$ be a solution to integral equation~(I) in ${\bf R}^N\times[0,1)$ 
with $\mu$ replaced by 
$$
\mu_T(x):=T^{\frac{p}{p-1}}\mu(T^{\frac{1}{\theta}}x),\qquad x\in{\bf R}^N.
$$ 
Then it follows from \eqref{eq:1.6} that  
$$
\sup_{x\in{\bf R}^N}\|\mu_T\|_{L^r(B(x,1))}\le \gamma.
$$
Set 
$$
u(x,t):=T^{-\frac{1}{p-1}}u_T(T^{-\frac{1}{\theta}}x,T^{-1}t),\qquad (x,t)\in{\bf R}^N\times(0,T).
$$
Then $u$ is a solution to integral equation~(I) in ${\bf R}^N\times[0,T)$. 
Therefore, for the proof of Theorem~\ref{Theorem:1.3}, 
we have only to prove the existence of solutions in ${\bf R}^N\times[0,1)$ 
under assumption~\eqref{eq:1.6} with $\sigma=1$.

For $1<r<\infty$, define
\[
\begin{aligned}
 & L^r_{uloc}:=\left\{ \phi\in L^r_{loc}({\bf R}^N)): \|\phi\|_{L^r_{uloc}}<\infty \right\},
 \quad
 \|\phi\|_{L^r_{uloc}}:= 
\sup_{x\in{\bf R}^N} \left( 
\int_{|x-y|<1} |\phi(y)|^r \,dy \right)^\frac{1}{r},\\
 & X_r:= L^\infty (0,1; L^r_{uloc}\,),\quad
\|f\|_{X_r}:=\underset{0<t<1}{\mbox{ess sup}}\,\|f(t)\|_{L^r_{uloc}}
\mbox{ for }f\in X_r.
\end{aligned}
\]
For any nonnegative Radon measure $\mu$ in ${\bf R}^N$, we set 
$$
[S(t)\mu](x):=\int_{{\bf R}^N}G(x-y,t)\,d\mu(y),\qquad x\in{\bf R}^N.
$$
Similarly, 
for any $f\in L^r_{uloc}$, we define
$$
[S(t)f](x):=\int_{{\bf R}^N}G(x-y,t)f(y)\,dy,\qquad x\in{\bf R}^N.
$$
Then we have (see \cite[Corollary 3.1]{MT06}):
\begin{lemma}
\label{Lemma:3.2}
Let $1\le r\le q\le \infty$. 
Then there exists $C=C(N,\theta)>0$ such that 
\begin{equation}
\label{eq:3.8}
\|S(t)f\|_{L^q_{uloc}} \le Ct^{-\frac{N}{\theta}\left(\frac{1}{r}-\frac{1}{q}\right)}
\|f\|_{L^r_{uloc}},
\qquad 0<t\le 1,  
\end{equation}
for $f\in L^r_{uloc}$.
\end{lemma}
Similarly, we have 
\begin{equation}
\label{eq:3.9}
\|S(t)\mu\|_{L^q_{uloc}} 
\le Ct^{-\frac{N}{\theta}\left(1-\frac{1}{q}\right)}
\sup_{x\in{\bf R}^N}\mu(B(x,1)),
\qquad 0<t\le 1,
\end{equation}
for nonnegative Radon measures $\mu$ in ${\bf R}^N$ and $1\le q\le\infty$. 
See also \cite[Lemma 2.1]{HI18} for $q=\infty$. 
\vspace{5pt}

We are ready to prove Theorem~\ref{Theorem:1.3}. 
\vspace{3pt}
\newline
{\bf Proof of Theorem~\ref{Theorem:1.3}.}  
It suffices to consider the case of $\sigma=1$. 
Let $0<\theta\le 2$ and $p>p_*$.  
Let $r_*<r<\infty$ and assume \eqref{eq:1.6} with $\sigma=1$. 

Let $\delta>0$ and $B_\delta:=\{f\in X_{pr}\,:\,\|f\|_{X_{pr}}<\delta\}$. 
Taking sufficiently small $\delta>0$ and $\gamma>0$, 
we show that 
\begin{equation}
\label{eq:3.10}
B_\delta\ni u\mapsto\Phi[u]\in B_\delta,
\end{equation}
where
\begin{equation}
\label{eq:3.11}
\Phi[u](t):= \int_0^t  S(t-s)\mu\,ds +\int_0^t S(t-s) |u(s)|^p\,ds
\quad\mbox{for}\quad t\in(0,1). 
\end{equation}
It follows from $r>r_*$ that 
$$
\frac{N}{\theta r}\left(1-\frac{1}{p}\right)<1,
$$
which together with Lemma~\ref{Lemma:3.2} and \eqref{eq:1.6} implies that
\begin{equation}
\label{eq:3.12}
\begin{split}
\left\|\,\int_0^t  S(t-s)\mu\,ds\,\right\|_{L^{pr}_{uloc}} 
 & \leq  \int_0^t \|  S(t-s)\mu \|_{L^{pr}_{uloc}}  \,ds\\
 & \le C\|\mu\|_{L^r_{uloc}}\int_0^t(t-s)^{-\frac{N}{\theta r}\left(1-\frac{1}{p}\right)}\,ds
 \le C\gamma 
\end{split}
\end{equation}
for $0<t<1$. 
Furthermore, 
by Lemma~\ref{Lemma:3.2}, for any $u\in B_\delta$, 
we have 
\begin{equation}
\label{eq:3.13}
\begin{split}
 & \left\|\,\int_0^t  S(t-s)|u(s)|^p\,ds\,\right\|_{L^{pr}_{uloc}}  
\le \int_0^t  \| S(t-s)|u(s)|^p\|_{L^{pr}_{uloc}}\,ds \\
 & \le C \int_0^t
 (t-s)^{-\frac{N}{\theta r}\left(1-\frac{1}{p}\right)}\|u(s)^p\|_{L^r_{uloc}}\,ds
=C \int_0^t (t-s)^{-\frac{N}{\theta r}\left(1-\frac{1}{p}\right)}\|u(s)\|_{L^{pr}_{uloc}}^p\,ds \\
 & \le C \delta^p \int_0^t  (t-s)^{-\frac{N}{\theta r}\left(1-\frac{1}{p}\right)}\,ds 
\le C \delta^p
\end{split}
\end{equation}
for $0<t<1$. 
By \eqref{eq:3.12} and \eqref{eq:3.13}, 
taking sufficiently small $\delta>0$ and $\gamma>0$ if necessary, 
we obtain 
$$
\|\Phi[u]\|_{X_{pr}}\le 
C\gamma+C \delta^p\le\delta,
$$
which implies \eqref{eq:3.10}.

Similarly, for any $u_1$, $u_2\in B_\delta$, 
we have 
\begin{equation}
\label{eq:3.14}
\begin{split}
\|\Phi[u_1](t)-\Phi[u_2](t)\|_{L^{pr}_{uloc}}
 & \le\int_0^t\|S(t-s)[|u_1(s)|^p-|u_2(s)|^p]\|_{L^{pr}_{uloc}}\,ds\\
 & \le C\int_0^t (t-s)^{-\frac{N}{\theta r}\left(1-\frac{1}{p}\right)}\||u_1(s)|^p-|u_2(s)|^p\|_{L^r_{uloc}}\,ds
\end{split}
\end{equation}
for $0<t<1$. 
Since 
$$
||u_1(x,t)|^p-|u_2(x,t)|^p|\le C(|u_1(x,t)|^{p-1}+|u_2(x,t)|^{p-1})|u_1(x,t)-u_2(x,t)|
$$
for $(x,t)\in{\bf R}^N\times(0,1)$, 
by the H\"older inequality we obtain 
\begin{equation*}
\begin{split}
\|u_1(s)^p-u_2(s)^p\|_{L^r_{uloc}} 
 & \le C \| u_1(s)-u_2(s) \|_{L^{pr}_{uloc}} 
\left( \| u_1(s) \|^{p-1}_{L^{pr}_{uloc}} + \| u_2(s) \|^{p-1}_{L^{pr}_{uloc}} \right) \\
 & \le C\delta^{p-1}\|u_1-u_2\|_X
\end{split}
\end{equation*}
for $0<s<1$.
This together with \eqref{eq:3.14} implies that 
\begin{equation}
\label{eq:3.15}
\|\Phi[u_1]-\Phi[u_2]\|_{X_{pr}}
\le C\delta^{p-1}\|u_1-u_2\|_{X_{pr}}.
\end{equation}
By \eqref{eq:3.10} and \eqref{eq:3.15}, 
taking a sufficiently small $\delta>0$ if necessary, 
we apply the contraction mapping theorem to obtain a fixed point $u\in B_\delta$. 
Then $u\ge 0$ in ${\bf R}^N\times(0,1)$ and it is a solution to integral equation~(I) in ${\bf R}^N\times[0,1)$. 
Thus Theorem~\ref{Theorem:1.3} follows.
$\Box$
\vspace{3pt}

Similarly to the proof of Theorem~\ref{Theorem:1.3}, 
we prove Theorem~\ref{Theorem:1.2}.
\vspace{3pt}
\newline
{\bf Proof of Theorem~\ref{Theorem:1.2}.} 
Let $1<p<p_*$. Then 
$$
r_*=\frac{N(p-1)}{\theta p}<1.
$$
We apply a similar argument to that of the proof of Theorem~\ref{Theorem:1.3} 
with $r$ replaced by $1$. 
It suffices to consider the case of $\sigma=1$. 
Using \eqref{eq:3.9}, instead of \eqref{eq:3.8}, 
and taking sufficiently small $\delta>0$ and $\gamma>0$ if necessary, 
we show that $\Phi$ defined by \eqref{eq:3.11} 
is a contraction mapping on $B_\delta\subset X_p$. 
Similarly to the proof of Theorem~\ref{Theorem:1.3}, 
we find a solution to integral equation~(I) in ${\bf R}^N\times[0,1)$. 
Thus Theorem~\ref{Theorem:1.2} follows. 
$\Box$
\section{Preliminary lemmas for the proof of Theorem~\ref{Theorem:1.4}}
In this section we prove some lemmas used in the proof of Theorem~\ref{Theorem:1.4}. 
For a nonnegative measurable function $\mu=\mu(x)$ in ${\bf R}^N$, 
set
\begin{equation}
\label{eq:4.1}
U(x,t)= \int_0^t \int_{{\bf R}^N} G(x-y,t-s) \mu(y) \,dy \,ds,
\end{equation}
which is a volume potential of $\mu$. 
We construct a solution to \eqref{eq:1.1} in the form $u=U+v$, 
that is, $v$ is a solution to 
\begin{equation}
\label{eq:4.2}
\left\{ 
\begin{aligned}
	&\partial_t v + (-\Delta)^\frac{\theta}{2} v= (U+v)^p 
	&&  \mbox{ in } {\bf R}^N\times (0,T), \\
	&v(\cdot,0)=0
	&&  \mbox{ in } {\bf R}^N.   
\end{aligned}
\right.
\end{equation}
The notion of solutions and supersolutions to problem~\eqref{eq:4.2} 
is defined in the same way as in Definition \ref{Definition:3.1}. 
Remark that the corresponding integral equation to \eqref{eq:4.2} is 
\begin{equation}
\label{eq:4.3}
\infty>v(x,t)=\Psi[v](x,t),
\end{equation}
where 
\begin{equation}
\label{eq:4.4}
\Psi[v](x,t) :=
\int_0^t \int_{{\bf R}^N} G(x-y,t-s) (U(y,s)+v(y,s))^p \,dy\,ds. 
\end{equation}
The proof of Theorem~\ref{Theorem:1.4} is based on the following lemma. 
\begin{lemma}\label{Lemma:4.1}
Assume that there exists a supersolution $\overline{v}$ 
to problem~\eqref{eq:4.2} in ${\bf R}^N\times[0,T)$, 
where $0<T\le\infty$. 
Then there exists a solution $v$ to problem~\eqref{eq:4.2} in ${\bf R}^N\times[0,T)$. 
In particular, $U+v$ is a solution to problem~\eqref{eq:1.1} in ${\bf R}^N\times[0,T)$. 
\end{lemma} 
{\bf Proof.}
Define a sequence $\{v_n\}_{n\geq 1}$ by 
$v_1:=\Psi[0]$ and $v_n:=\Psi[v_{n-1}]$. 
Since $U\geq0$ and $\Psi[\overline{v}]\leq \overline{v}$, 
by induction we see that 
$$
0\leq v_1(x,t)\leq v_2(x,t) \leq \cdots \leq v_n(x,t) \leq \cdots \leq \overline{v}(x,t)< \infty
$$
for almost all $(x,t)\in{\bf R}^N\times(0,T)$. Then
$v(x,t):=\lim_{n\to\infty} v_n(x,t)$ is well-defined for almost all $(x,t)\in {\bf R}^N\times(0,T)$ 
and $v$ satisfies \eqref{eq:4.3}, that is, $v$ is a solution to problem~\eqref{eq:4.2} in ${\bf R}^N\times[0,T)$. 
Then $U+v$ is a solution to problem~\eqref{eq:1.1} in ${\bf R}^N\times[0,T)$. 
Thus Lemma~\ref{Lemma:4.1} follows. 
$\Box$
\vspace{5pt}

In the rest of this section we obtain two lemmas on volume potentials and recall some properties of $G$. 
By the letter $C$
we denote generic positive constants independent of $x\in{\bf R}^N$ and $t>0$ 
and they may have different values also within the same line. 
\begin{lemma}
\label{Lemma:4.2}
Let $N\geq1$ and $0<\theta\leq 2$. 
Let $a\geq0$ and $b<N$ satisfy $(1+a)\theta<b$. 
Then there exists a constant $C=C(N,\theta,a,b)>0$ such that 
\begin{equation}
\label{eq:4.5}
\int_0^t \int_{{\bf R}^N} G(x-y,t-s) (t-s)^a |y|^{-b} \,dy\,ds 
\leq C|x|^{(1+a)\theta - b }
\end{equation}
for $x\in{\bf R}^N\setminus\{0\}$ and $t>0$. 
\end{lemma}
{\bf Proof.}
It follows from \eqref{eq:3.2} and \eqref{eq:3.4} that 
$$
0\le G(x,t)\le Ct^{-\frac{N}{\theta}} 
(1+t^{-\frac{1}{\theta}}|x|)^{-N-\theta}
$$
for $(x,t)\in{\bf R}^N\times(0,\infty)$,
where $0<\theta\le 2$. 
Since $(1+a)\theta<b<N$, 
we have
\begin{equation}
\label{eq:4.6}
\begin{split}
 & \int_0^t G(x-y,t-s) (t-s)^a\,ds\\
 & \le C\int_0^t (t-s)^{-\frac{N}{\theta}+a}(1+(t-s)^{-\frac{1}{\theta}}|x-y|)^{-N-\theta}\,ds\\
 & \le C|x-y|^{(1+a)\theta-N}\int_0^\infty \eta^{\frac{N}{\theta}-a-2} 
	( 1+\eta^\frac{1}{\theta} )^{-N-\theta}\,d\eta\le C|x-y|^{(1+a)\theta-N}	
\end{split}
\end{equation}
for $x$, $y\in{\bf R}^N$ with $x\not=y$ and $t>0$. 
Then it follows that
\begin{equation}
\label{eq:4.7}
\int_0^t \int_{{\bf R}^N} G(x-y,t-s) (t-s)^a |y|^{-b} \,dy\,ds
\le C\sum_{i=1}^4\int_{\Omega_i}|x-y|^{(1+a)\theta-N}|y|^{-b}\,dy
\end{equation}
for $x\in{\bf R}^N\setminus\{0\}$ and $t>0$, where 
\begin{equation}
\label{eq:4.8}
\begin{aligned}
	& \Omega_1:= \{|y|\leq |x|/2\}, 
	&& \Omega_2:=\{|x-y|\leq |x|/2\}, \\
	& \Omega_3:=\{|y|\geq |x|/2, |y|\leq |x-y|\}, 
	&& \Omega_4:=\{|x-y|\geq |x|/2, |y|\geq |x-y|\}. 
\end{aligned}
\end{equation}
Then we have: 
\begin{itemize}
  \item[($\Omega_1$)] 
  Let $y\in\Omega_1$. Then $|x-y|\ge |x|-|y|\ge |x|/2$. 
  This implies that 
  $$
  \int_{\Omega_1}|x-y|^{(1+a)\theta-N}|y|^{-b}\,dy
  \le C|x|^{(1+a)\theta-N}
  \int_{\Omega_1}  |y|^{-b}\,dy
  \le C|x|^{(1+a)\theta-b};
  $$
  \item[($\Omega_2$)] 
  Let $y\in\Omega_2$. Then $|y|\ge |x|-|y-x|\ge |x|/2$ and 
  $$
  \int_{\Omega_2}|x-y|^{(1+a)\theta-N}|y|^{-b}\,dy\le C|x|^{-b}\int_{\Omega_2}|x-y|^{(1+a)\theta-N}\,dy
  \le C|x|^{(1+a)\theta-b};
  $$
  \item[($\Omega_3$)] 
  Let $y\in\Omega_3$. Then 
  $$
  \int_{\Omega_3}|x-y|^{(1+a)\theta-N}|y|^{-b}\,dy\le C\int_{\{|y|\ge |x|/2\}}|y|^{(1+a)\theta-b-N}\,dy
  \le C|x|^{(1+a)\theta-b};
  $$
  \item[($\Omega_4$)] 
  Let $y\in\Omega_4$. Then
  $$
  \int_{\Omega_4}|x-y|^{(1+a)\theta-N}|y|^{-b}\,dy
  \le\int_{\{|x-y|\ge|x|/2\}}|x-y|^{(1+a)\theta-b-N}\,dy
  \le C|x|^{(1+a)\theta-b}.
  $$
\end{itemize}
These together with \eqref{eq:4.7} imply \eqref{eq:4.5}. 
Thus Lemma~\ref{Lemma:4.2} follows. 
$\Box$
\begin{lemma}
\label{Lemma:4.3}
Let $N\geq1$ and $0<\theta\leq 2$ satisfy $\theta<N$. 
Let $a\geq0$ and $b>1$. 
Then there exist $C=C(N,\theta,a,b)>0$ and $A>1$ such that 
\begin{equation}
\label{eq:4.9}
\begin{split}
 & \int_0^t \int_{{\bf R}^N} G(x-y,t-s) \left[\log\biggr(A+\frac{1}{t-s}\biggr) \right]^{-a} 
 |y|^{-N}\left[\log\biggr(A+\frac{1}{|y|}\biggr)\right]^{-b} \,dy\,ds \\
 & \le C|x|^{\theta-N} \left[\log\biggr(A+\frac{1}{|x|}\biggr)\right]^{-a-b+1}\log(A+t)
\end{split}
\end{equation}
for $x\in{\bf R}^N\setminus\{0\}$ and $t>0$. 
\end{lemma}
{\bf Proof.}
The proof is divided into 3 steps. 
Let $A>1$ be sufficiently large. 
\newline
\underline{Step 1.} 
Let $0<\delta<2\theta$. We prove that 
\begin{equation}
\label{eq:4.10}
\begin{split}
 & \int_0^t G(x-y,t-s) \left[\log\biggr(A+\frac{1}{t-s}\biggr) \right]^{-a}\,ds\\
 & \le C|x-y|^{\theta-N} \left[\log\biggr(A+\frac{1}{|x-y|^\theta}\biggr) \right]^{-a}
 (1+t^{-\frac{1}{\theta}}|x-y| )^{-\delta}
\end{split}
\end{equation}
for $x$, $y\in{\bf R}^N$ with $x\not=y$ and $t>0$. 
To this end, similarly to \eqref{eq:4.6}, we see that 
\begin{equation}
\label{eq:4.11}
\begin{split}
 & \int_0^t G(x-y,t-s) \left[\log\biggr(A+\frac{1}{t-s}\biggr) \right]^{-a}\,ds\\
 & \le C\int_0^t (t-s)^{-\frac{N}{\theta}}\left[\log\biggr(A+\frac{1}{t-s}\biggr) \right]^{-a}
 (1+(t-s)^{-\frac{1}{\theta}}|x-y|)^{-N-\theta}\,ds\\
 & \le C|x-y|^{\theta-N} 
\int_{t^{-1}|x-y|^\theta}^\infty \eta^{\frac{N}{\theta}-2} 
(1 + \eta^\frac{1}{\theta})^{-N-\theta}
\left[\log\biggr(A+\frac{\eta}{|x-y|^\theta}\biggr) \right]^{-a}\,d\eta\\
 & \le C|x-y|^{\theta-N} (1+t^{-\frac{1}{\theta}}|x-y| )^{-\delta}\\
 & \qquad\quad
 \times\int_0^\infty \eta^{\frac{N}{\theta}-2} 
(1 + \eta^\frac{1}{\theta})^{-N-\theta+\delta}
\left[\log\biggr(A+\frac{\eta}{|x-y|^\theta}\biggr) \right]^{-a}\,d\eta
\end{split}
\end{equation}
for $x$, $y\in{\bf R}^N$ with $x\not=y$ and $t>0$. 
On the other hand, 
\begin{equation}
\label{eq:4.12}
\begin{split}
 & \int_1^\infty \eta^{\frac{N}{\theta}-2} 
(1 + \eta^\frac{1}{\theta})^{-N-\theta+\delta}
\left[\log\biggr(A+\frac{\eta}{|x-y|^\theta}\biggr) \right]^{-a}\,d\eta\\
 & \le\left[\log\biggr(A+\frac{1}{|x-y|^\theta}\biggr) \right]^{-a}
\int_1^\infty \eta^{\frac{N}{\theta}-2} 
(1 + \eta^\frac{1}{\theta})^{-N-\theta+\delta}\,d\eta\\
 & \le C\left[\log\biggr(A+\frac{1}{|x-y|^\theta}\biggr) \right]^{-a}
\end{split}
\end{equation}
for $x$, $y\in{\bf R}^N$ with $x\not=y$.
Let $\epsilon>0$ be sufficiently small such that $0<\epsilon<(N/\theta)-1$. 
Taking a sufficiently large $A>1$ if necessary, we see that 
$r\mapsto r^\epsilon[\log(A +r)]^{-a}$ is increasing. 
Then 
\begin{equation}
\label{eq:4.13}
\begin{split}
 & \int_0^1\eta^{\frac{N}{\theta}-2}
 (1 + \eta^\frac{1}{\theta})^{-N-\theta+\delta}
 \left[\log\biggr(A+\frac{\eta}{|x-y|^\theta}\biggr) \right]^{-a}\,d\eta\\
 & \le C\int_0^1\eta^{\frac{N}{\theta}-2}
 \left(\frac{\eta}{|x-y|^\theta}\right)^{-\epsilon}
 \left(\frac{\eta}{|x-y|^\theta}\right)^\epsilon
 \left[\log\biggr(A+\frac{\eta}{|x-y|^\theta}\biggr) \right]^{-a}\,d\eta\\
 & \le C\int_0^1\eta^{\frac{N}{\theta}-2}
 \left(\frac{\eta}{|x-y|^\theta}\right)^{-\epsilon}
 \left(\frac{1}{|x-y|^\theta}\right)^\epsilon
 \left[\log\biggr(A+\frac{1}{|x-y|^\theta}\biggr) \right]^{-a}\,d\eta\\
 & \le C\left[\log\biggr(A+\frac{1}{|x-y|^\theta}\biggr) \right]^{-a}
 \int_0^1\eta^{\frac{N}{\theta}-2-\epsilon}\,d\eta
 \le C\left[\log\biggr(A+\frac{1}{|x-y|^\theta}\biggr) \right]^{-a}
\end{split}
\end{equation}
for $x$, $y\in{\bf R}^N$ with $x\not=y$.
Combining \eqref{eq:4.11}, \eqref{eq:4.12} and \eqref{eq:4.13}, 
we obtain \eqref{eq:4.10}. 
\vspace{3pt}
\newline
\underline{Step 2.} 
It follows from \eqref{eq:4.10} that 
\begin{equation}
\label{eq:4.14}
\begin{split}
 & \int_0^t \int_{{\bf R}^N} G(x-y,t-s) \left[\log\biggr(A+\frac{1}{t-s}\biggr) \right]^{-a} 
 |y|^{-N}\left[\log\biggr(A+\frac{1}{|y|}\biggr)\right]^{-b} \,dy\,ds \\
 & \le C\int_{{\bf R}^N}K(x,y,t)\,dy,
\end{split}
\end{equation}
where
\begin{equation}
\label{eq:4.15}
\begin{split}
 & K(x,y,t):=|x-y|^{\theta-N} \left[\log\biggr(A+\frac{1}{|x-y|^\theta}\biggr) \right]^{-a}\\
 & \qquad\qquad\qquad
\times (1+t^{-\frac{1}{\theta}}|x-y| )^{-\delta}
|y|^{-N}\left[\log\biggr(A+\frac{1}{|y|}\biggr)\right]^{-b}.
\end{split}
\end{equation}
Let $x\in B(0,1)\setminus\{0\}$ and $t>0$. We prove that
\begin{equation}
\label{eq:4.16}
\int_{{\bf R}^N}K(x,y,t)\,dy
=\sum_{i=1}^4\int_{\Omega_i}K(x,y,t)\,dy\le C|x|^{\theta-N}
\left[\log\biggr(A+\frac{1}{|x|}\biggr)\right]^{-a-b+1},
\end{equation}
where $\Omega_i$ $(i=1,2,3,4)$ is as in \eqref{eq:4.8}. 
Taking a sufficiently large $A>1$ if necessary, 
we see that 
$r\mapsto r^{\theta-N}[\log(A +r^{-\theta})]^{-a}$ is decreasing. 
Then we have 
\begin{equation*}
\begin{split}
 & \int_{\Omega_1}K(x,y,t)\,dy\\
 & \le C\int_{\Omega_1}\left(\frac{|x|}{2}\right)^{\theta-N} \left[\log\biggr(A+\frac{2^\theta}{|x|^\theta}\biggr) \right]^{-a}
|y|^{-N}\left[\log\biggr(A+\frac{1}{|y|}\biggr)\right]^{-b}\,dy\\
 & \le C|x|^{\theta-N} \left[\log\biggr(A+\frac{1}{|x|}\biggr) \right]^{-a}
 \int_0^{|x|/2} r^{-1}\left[\log\biggr(A+\frac{1}{r}\biggr)\right]^{-b}\,dr\\
 & \le C|x|^{\theta-N} \left[\log\biggr(A+\frac{1}{|x|}\biggr) \right]^{-a-b+1}. 
\end{split}
\end{equation*}
(See also $(\Omega_1)$ in the proof of Lemma~\ref{Lemma:4.2}.)
Similarly, taking a sufficiently large $A>1$ if necessary,  
we obtain the following estimates 
(see also $(\Omega_2)$, $(\Omega_3)$ and $(\Omega_4)$ in the proof of Lemma~\ref{Lemma:4.2}): 
\begin{equation*}
\begin{split}
\,\,\, & \int_{\Omega_2}K(x,y,t)\,dy\\
 & \le C\int_{\Omega_2}|x-y|^{\theta-N} \left[\log\biggr(A+\frac{1}{|x-y|^\theta}\biggr) \right]^{-a}
|x|^{-N}\left[\log\biggr(A+\frac{1}{|x|}\biggr)\right]^{-b}\,dy\\
 & \le C|x|^{-N}\left[\log\biggr(A+\frac{1}{|x|}\biggr)\right]^{-b}
 \int_0^{|x|/2}r^{\theta-1} \left[\log\biggr(A+\frac{1}{r^\theta}\biggr) \right]^{-a}\,dr\\
 & \le C|x|^{\theta-N}\left[\log\biggr(A+\frac{1}{|x|}\biggr)\right]^{-a-b};
\end{split}
\end{equation*}
\begin{equation*}
\begin{split}
 & \int_{\Omega_3}K(x,y,t)\,dy\\
 & \le C\int_{\Omega_3}|y|^{\theta-N} \left[\log\biggr(A+\frac{1}{|y|^\theta}\biggr) \right]^{-a}
|y|^{-N}\left[\log\biggr(A+\frac{1}{|y|}\biggr)\right]^{-b}\,dy\qquad\quad\\
 & \le C|x|^{\theta-N+\epsilon'}
 \left[\log\biggr(A+\frac{1}{|x|}\biggr)\right]^{-a-b}
 \int_{|x|/2}^\infty r^{-1-\epsilon'}\,dr\\
 & \le C|x|^{\theta-N}\left[\log\biggr(A+\frac{1}{|x|}\biggr)\right]^{-a-b};
\end{split}
\end{equation*}
\begin{equation*}
\begin{split}
 & \int_{\Omega_4}K(x,y,t)\,dy\\
 & \le C\int_{\Omega_4}|x|^{\theta-N} \left[\log\biggr(A+\frac{2^\theta}{|x|^\theta}\biggr) \right]^{-a}
|x-y|^{-N}\left[\log\biggr(A+\frac{1}{|x-y|}\biggr)\right]^{-b}\,dy\\
 & \le C|x|^{\theta-N+\epsilon'} \left[\log\biggr(A+\frac{1}{|x|}\biggr) \right]^{-a-b}
 \int_{|x|/2}^\infty r^{-1-\epsilon'} \,dr\\
 & \le C|x|^{\theta-N}\left[\log\biggr(A+\frac{1}{|x|}\biggr)\right]^{-a-b}.
\end{split}
\end{equation*}
Here $0<\epsilon'<N-\theta$. 
These imply \eqref{eq:4.16}.
\vspace{3pt}
\newline
\underline{Step 3.} 
Let $x\in{\bf R}^N\setminus B(0,1)$ and $t>0$. 
We prove that
\begin{equation}
\label{eq:4.17}
\int_{{\bf R}^N}K(x,y,t)\,dy
=\sum_{i=1}^4\int_{\tilde{\Omega}_i}K(x,y,t)\,dy
\le C|x|^{\theta-N}\log(A+t),
\end{equation}
where 
\[
\begin{aligned}
	&\tilde \Omega_1:= \{|y|\leq 1/2\}, 
	&&\tilde \Omega_2:=\{|x-y|\leq 1/2 \}, \\
	&\tilde \Omega_3:=\{|y|\geq 1/2, |y|\leq |x-y|\}, 
	&&\tilde \Omega_4:=\{|x-y|\geq 1/2, |y|\geq |x-y|\}.
\end{aligned}
\]
Since $|x-y|\ge |x|/2\ge 1/2$ for $y\in\tilde{\Omega}_1$, 
we observe from \eqref{eq:4.15} and $b>1$ that 
\begin{equation}
\label{eq:4.18}
\int_{\tilde{\Omega}_1}K(x,y,t)\,dy
\le C|x|^{\theta-N}\int_{\tilde{\Omega}_1}
|y|^{-N}\left[\log\biggr(A+\frac{1}{|y|}\biggr)\right]^{-b}\,dy
\le C|x|^{\theta-N}.
\end{equation}
On the other hand, 
it follows that $|y|\ge|x|-|y-x|\ge |x|/2$ and $|x-y|\leq |x|+|y|\leq 3|y|$ for $y\in \tilde\Omega_2$. 
Taking a sufficiently large $A>1$ if necessary, 
we see that
$s\mapsto s^{-\theta}[\log(A+s^{-1})]^{-b}$ is decreasing. 
Then we have 
\begin{equation}
\label{eq:4.19}
\begin{split}
 & \int_{\tilde{\Omega}_2}K(x,y,t)\,dy\\
 & \le C|x|^{\theta-N}\int_{\tilde{\Omega}_2}|x-y|^{\theta-N}\left[\log\biggr(A+\frac{1}{|x-y|^\theta}\biggr) \right]^{-a}
|y|^{-\theta}\left[\log\biggr(A+\frac{1}{|y|}\biggr)\right]^{-b}\,dy\\
 & \le C|x|^{\theta-N}
\int_{\{|x-y|\le 1/2\}}|x-y|^{-N}\left[\log\biggr(A+\frac{1}{|x-y|}\biggr)\right]^{-a-b}\,dy
\le C|x|^{\theta-N}.
\end{split}
\end{equation}
For $y\in\tilde{\Omega}_3$, 
it follows that $|x-y|\geq |x|-|y|\geq |x|-|x-y|$, 
that is,  $|x-y|\geq |x|/2\ge 1/2$. 
Then we obtain
\begin{equation}
\label{eq:4.20}
\begin{split}
\int_{\tilde{\Omega}_3}K(x,y,t)\,dy
 & \le C|x|^{\theta-N}\int_{\tilde{\Omega}_3}
|y|^{-N} (1+t^{-\frac{1}{\theta}}|y| )^{-\delta}\,dy\\
 & \le C|x|^{\theta-N}\int_{1/2}^\infty
r^{-1} (1+t^{-\frac{1}{\theta}}r )^{-\delta}\,dr
\le C|x|^{\theta-N}\log(A+t).
\end{split}
\end{equation}
Furthermore, for $y\in \tilde\Omega_4$, we have $|y|\ge|x-y|\ge|x|-|y|$, that is, $|y|\ge |x|/2$. 
Then we observe that 
\begin{equation}
\label{eq:4.21}
\begin{split}
\int_{\tilde{\Omega}_4}K(x,y,t)\,dy
 & \le C|x|^{\theta-N}\int_{\tilde{\Omega}_4}
|x-y|^{\theta-N}|y|^{-\theta} (1+t^{-\frac{1}{\theta}}|y| )^{-\delta}\,dy\\
 & \le C|x|^{\theta-N}\int_{\tilde{\Omega}_4}
|x-y|^{\theta-N}|x-y|^{-\theta} (1+t^{-\frac{1}{\theta}}|x-y| )^{-\delta}\,dy\\
 & \le C|x|^{\theta-N}\int_{1/2}^\infty r^{-1} (1+t^{-\frac{1}{\theta}}r )^{-\delta}\,dr
\le C|x|^{\theta-N}\log(A+t).
\end{split}
\end{equation}
Combining \eqref{eq:4.18}, \eqref{eq:4.19}, \eqref{eq:4.20} and \eqref{eq:4.21}, 
we obtain \eqref{eq:4.17}. 
Therefore, 
by \eqref{eq:4.14}, \eqref{eq:4.16} and \eqref{eq:4.17} 
we have \eqref{eq:4.9}. Thus Lemma~\ref{Lemma:4.3} follows.
$\Box$
\vspace{5pt}
\newline
Furthermore, we recall two lemmas on the fundamental solution~$G$.
\begin{lemma}
\label{Lemma:4.4}
There exists $C=C(N,\theta)>0$ such that 
\[
\begin{aligned}
	& \frac{G(x-y,t-s) G(y-z,s-\tau)}{G(x-z,t-\tau)}\\
	& \left\{ 
	\begin{aligned}
	&=
	G\left( y-\frac{s-\tau}{t-\tau}x - \frac{t-s}{t-\tau}z, 
	\frac{(t-s)(s-\tau)}{t-\tau}\right) 
	&& \mbox{ if }\quad\theta=2, \\
	&\leq 
	C \left( G(x-y,t-s) + G(y-z,s-\tau) \right) 
	&& \mbox{ if }\quad0<\theta<2, 
	\end{aligned}
	\right.
\end{aligned}
\]
for $x$, $y$, $z\in {\bf R}^N$ and $\tau<s<t$. 
\end{lemma}
{\bf Proof.}
In the case of $\theta=2$, Lemma~\ref{Lemma:4.4} is shown by straightforward computations. 
In the case of $0<\theta<2$, Lemma~\ref{Lemma:4.4} follows from \cite[Theorem~4]{BJ07}. 
$\Box$

\begin{lemma}
\label{Lemma:4.5}
Let $0<a<N$ and $b\ge 0$. 
Then there exist $C=C(N,\theta,a,b)>0$ and $A>1$ such that 
\begin{equation}
\label{eq:4.22}
\sup_{x\in{\bf R}^N}\int_{{\bf R}^N}G(x-y,t)
|y|^{-a}\left[\log\biggr(A+\frac{1}{|y|}\biggr) \right]^{-b}
\,dy\le Ct^{-\frac{a}{\theta}}\left[\log\biggr(A+\frac{1}{t}\biggr) \right]^{-b}
\end{equation}
for $t>0$. 
\end{lemma}
{\bf Proof.}
Let $\delta>0$ satisfy $a+\delta<N$. 
We choose $A>1$ so large that 
$\phi(s):=s^{-a}[\log(A+s^{-1})]^{-b}$ is decreasing 
and $s\mapsto s^\delta[\log(A+s^{-1}) ]^{-b}$ is increasing. 

It follows from \cite[Lemma~2.1]{HI18} that
\begin{equation*}
\sup_{x\in{\bf R}^N}
\int_{{\bf R}^N}G(x-y,t)\phi(|y|)\,dy\le Ct^{-\frac{N}{\theta}}
\sup_{x\in{\bf R}^N}\int_{B(x,t^{1/\theta})}\phi(|y|)\,dy
\end{equation*}
for $t>0$. For $x\in B(0,2t^{1/\theta})$, we have 
\begin{equation*}
\begin{split}
 & \int_{B(x,t^{1/\theta})}\phi(|y|)\,dy
\le\int_{B(0,3t^{1/\theta})}\phi(|y|)\,dy\\
 &\quad  \le (3t^{\frac{1}{\theta}})^\delta\left[\log\biggr(A+\frac{1}{3t^{\frac{1}{\theta}}}\biggr) \right]^{-b}
\int_{B(0,3t^{1/\theta})}|y|^{-a-\delta}\,dy
 \le Ct^{\frac{N}{\theta}-\frac{a}{\theta}}
 \left[\log\biggr(A+\frac{1}{t}\biggr) \right]^{-b} 
\end{split}
\end{equation*}
for $t>0$. 
On the other hand, 
since $|y|\geq |x|-|x-y|\geq t^{1/\theta}$ for $x\in {\bf R}^N\setminus B(0,2t^{1/\theta})$ 
and $y\in B(x,t^{1/\theta})$, we obtain
\begin{equation*}
\int_{B(x,t^{1/\theta})}\phi(|y|)\,dy
\le Ct^{\frac{N}{\theta}}\phi(t^{\frac{1}{\theta}})\le Ct^{\frac{N}{\theta}-\frac{a}{\theta}}\left[\log\biggr(A+\frac{1}{t}\biggr) \right]^{-b}
\end{equation*}
for $x\in {\bf R}^N\setminus B(0,2t^{1/\theta})$ and $t>0$. 
Then \eqref{eq:4.22} follows, and the proof is complete.
$\Box$
\section{Proof of Theorem~\ref{Theorem:1.4} with $\boldsymbol{p>p_*}$}
Let $p>p_*$ and $\alpha>1$. Set 
\begin{equation}
\label{eq:5.1}
W(x,t):=
\int_0^t \int_{{\bf R}^N} G(x-y,t-s) |y|^{-\frac{\theta p}{p-1}\alpha} 
(t-s)^{\frac{p}{p-1}(\alpha-1)} \,dy\,ds.
\end{equation}
The main purpose of this section is to prove the following proposition,
which concerns with the existence of supersolutions to integral equation~(I). 
\begin{proposition}
\label{Proposition:5.1}
Assume the same conditions as in Theorem~{\rm\ref{Theorem:1.4}}.
Let $p>p_*$ and $W$ be as in \eqref{eq:5.1}. 
Then there exist $\gamma_*=\gamma_*(N,\theta,p)>0$ and $\alpha=\alpha(N,\theta,p)>1$ 
such that, for any $C_0\ge 0$, 
the function $w$ defined by 
$$
w(x,t):=\gamma_*W(x,t)+C_0t
$$ 
is a supersolution to integral equation~{\rm (I)} in ${\bf R}^N\times[0,T)$ for some $T\in(0,\infty]$. 
Furthermore, $T=\infty$ if $C_0=0$. 
\end{proposition}
In what follows, we set 
$$
V[f](x,t):=\int_0^t \int_{{\bf R}^N} G(x-y,t-s) f(y,s)\,dy\,ds
$$
for nonnegative measurable functions $f$ in ${\bf R}^N\times(0,\infty)$.
\vspace{3pt}
\newline
{\bf Proof.}
Let $p>p_*$ and assume \eqref{eq:1.7}. Let $x\in{\bf R}^N\setminus\{0\}$, $t>0$ and $U$ be as in \eqref{eq:4.1}. 
By the letter $C$
we denote generic positive constants which is independent of $x$ and $t$ 
and they may have different values also within the same line. 
The proof is divided into 3 steps. 
\vspace{3pt}
\newline
\underline{Step 1.} 
It follows from $p>p_*$ that $\theta p/(p-1)<N$. 
Let $\alpha>1$ be such that $\theta p\alpha/(p-1)<N$. 
We prove that  
\begin{equation}
\label{eq:5.2}
V[U^p](x,t)
\leq C\gamma^p W(x,t)+CC_0^pt^{p+1}.
\end{equation}
It follows from \eqref{eq:1.7}, \eqref{eq:3.5} and \eqref{eq:4.1} that 
\begin{equation}
\label{eq:5.3}
U(x,t)\le \gamma\int_0^t\int_{{\bf R}^N}G(x-y,t-s)|y|^{-\frac{\theta p}{p-1}}\,dy\,ds+C_0t.
\end{equation}
Let $1<k<p$.  
By the H\"older inequality and Lemma~\ref{Lemma:4.2}
we have 
\begin{equation*}
\begin{split}
 & \int_0^t\int_{{\bf R}^N}G(x-y,t-s)|y|^{-\frac{\theta p}{p-1}}\,dy\,ds\\
 & =\int_0^t\int_{{\bf R}^N}G(x-y,t-s)|y|^{-\frac{\theta p}{p-1}\frac{\alpha}{k}}
|y|^{-\frac{\theta p}{p-1}(1-\frac{\alpha}{k})}\,dy\,ds\\
 & \le\left(\int_0^t\int_{{\bf R}^N}G(x-y,t-s)|y|^{-\frac{\theta p}{p-1}\alpha}\,ds\,ds\right)^{\frac{1}{k}}\\
 & \qquad\quad
 \times\left(\int_0^t\int_{{\bf R}^N}G(x-y,t-s)|y|^{-\frac{\theta p}{p-1}\frac{k-\alpha}{k-1}}\,dy\,ds\right)^{1-\frac{1}{k}}\\
 & \le C\left(|x|^{\theta-\frac{\theta p}{p-1}\alpha}\right)^{\frac{1}{k}-\frac{1}{p}}
 \left(\int_0^t\int_{{\bf R}^N}G(x-y,t-s)|y|^{-\frac{\theta p}{p-1}\alpha}\,ds\,ds\right)^{\frac{1}{p}}
 \left(|x|^{\theta-\frac{\theta p}{p-1}\frac{k-\alpha}{k-1}}\right)^{1-\frac{1}{k}}\\
 & =C|x|^{-\frac{\lambda_1}{p}}\left(\int_0^t\int_{{\bf R}^N}G(x-y,t-s)|y|^{-\lambda_2}\,ds\,ds\right)^{\frac{1}{p}},
\end{split}
\end{equation*}
where 
$$
\lambda_1:= -\theta(p-1)+\frac{\theta p}{p-1}(p-\alpha)=\theta-\frac{\theta p}{p-1}(\alpha-1),\qquad 
\lambda_2:= \frac{\theta p}{p-1}\alpha. 
$$
We remark that $\lambda_1<\theta<N$ and $\lambda_2<N$.
Taking a suitable $\alpha>1$ if necessary, 
we can assume, without loss of generality, that 
\begin{equation}
\label{eq:5.4}
\lambda_1>0. 
\end{equation}
This together with \eqref{eq:5.3} implies that 
$$
U(x,t)^p\le C(C_0t)^p+C\gamma^p|x|^{-\lambda_1}
\int_0^t\int_{{\bf R}^N}G(x-y,t-s)|y|^{-\lambda_2}\,dy\,ds.
$$
Then we deduce from the Fubini theorem and \eqref{eq:3.5} that
\begin{equation}
\label{eq:5.5}
V[U^p](x,t)\leq C\gamma^p I(x,t)+CC_0^pt^{p+1},
\end{equation}
where
\begin{equation}
\label{eq:5.6}
\begin{split}
 & I(x,t)\\
 & :=\int_0^t \int_{{\bf R}^N}  |z|^{-\lambda_2}
 \left(\int_\tau^t \int_{{\bf R}^N} G(x-y,t-s) G(y-z,s-\tau)|y|^{-\lambda_1} \,dy\,ds\right)\,dz\,d\tau\\
 & \,=\int_0^t \int_{{\bf R}^N}  G(x-z,t-\tau)|z|^{-\lambda_2}J(x,z,t,\tau)\,dz\,d\tau
\end{split}
\end{equation} 
and 
$$
J(x,z,t,\tau):=
 \int_\tau^t \int_{{\bf R}^N}\frac{G(x-y,t-s) G(y-z,s-\tau)}{G(x-z,t-\tau)}|y|^{-\lambda_1} \,dy\,ds.
$$
In the case of $\theta=2$, 
by Lemmas~\ref{Lemma:4.4} and \ref{Lemma:4.5} 
we have
\begin{equation}
\label{eq:5.7}
\begin{split}
 & J(x,z,t,\tau)\\
 & \le C\int_\tau^t \int_{{\bf R}^N}G\left( y-\frac{s-\tau}{t-\tau}x - \frac{t-s}{t-\tau}z, 
\frac{(t-s)(s-\tau)}{t-\tau}\right)|y|^{-\lambda_1} \,dy\,ds\\
 & \le C\int_\tau^t \left\|\int_{{\bf R}^N}G(\cdot-y,\xi(s:t,\tau))|y|^{-\lambda_1} \,dy\right\|_{L^\infty({\bf R}^N)} \,ds
 \le C\int_\tau^t\xi(s:t,\tau)^{-\frac{\lambda_1}{2}}\,ds\\
 & \le C(t-\tau)^{\frac{\lambda_1}{2}}\int_\tau^t(t-s)^{-\frac{\lambda_1}{2}}(s-\tau)^{-\frac{\lambda_1}{2}}\,ds\\
 & \le C(t-\tau)^{-\frac{\lambda_1}{2}+1}=C(t-\tau)^{\frac{p}{p-1}(\alpha-1)},
\end{split}
\end{equation}
where
\begin{equation}
\label{eq:5.8}
\xi(s:t,\tau):=\frac{(t-s)(s-\tau)}{t-\tau}.
\end{equation}
Similarly, in the case of $0<\theta<2$, we have 
\begin{equation}
\label{eq:5.9}
\begin{split}
J(x,z,t,\tau)
 & \le C\int_\tau^t \int_{{\bf R}^N}[G(x-y,t-s) + G(y-z,s-\tau)]|y|^{-\lambda_1} \,dy\,ds\\
 & \le C\int_\tau^t (t-s)^{-\frac{\lambda_1}{\theta}}\,ds
 =C(t-\tau)^{-\frac{\lambda_1}{\theta}+1}=C(t-\tau)^{\frac{p}{p-1}(\alpha-1)}.
\end{split}
\end{equation}
Combining \eqref{eq:5.6}, \eqref{eq:5.7} and \eqref{eq:5.9}, we obtain 
\begin{equation}
\label{eq:5.10}
I(x,t)\le C\int_0^t \int_{{\bf R}^N}  G(x-z,t-\tau)|z|^{-\lambda_2}
(t-\tau)^{\frac{p}{p-1}(\alpha-1)}\,dz\,d\tau
=CW(x,t).
\end{equation}
This together with \eqref{eq:5.5} implies \eqref{eq:5.2}. 
\vspace{3pt}
\newline
\underline{Step 2.}
We show that 
\begin{equation}
\label{eq:5.11}
V[W^p](x,t)\le CW(x,t). 
\end{equation}
For this aim, 
by the H\"older inequality and Lemma~\ref{Lemma:4.2} 
we obtain 
\[
\begin{aligned}
	W(x,t)^p &\leq 
	\left(\int_0^t \int_{{\bf R}^N}
	G(x-y,t-s) |y|^{-\frac{\theta p}{p-1}\alpha} 
	(t-s)^{\left(\frac{p}{p-1}\right)^2 (\alpha-1)} \,dy\,ds  
	\right)^{p-1} \\
	&\quad \times 
	\int_0^t \int_{{\bf R}^N} G(x-y,t-s) |y|^{-\frac{\theta p}{p-1}\alpha } \,dy\,ds \\
	&\leq 
	C|x|^{-\lambda_1}
	\int_0^t \int_{{\bf R}^N} G(x-y,t-s) |y|^{-\lambda_2} \,dy\,ds. 
\end{aligned}
\]
Here we used $\lambda_1>0$ (see \eqref{eq:5.4}).
Then we deduce from \eqref{eq:5.6} and \eqref{eq:5.10} that
$$
V[W^p](x,t)\le CI(x,t)\le CW(x,t), 
$$
which implies \eqref{eq:5.11}. 
\vspace{3pt}
\newline
\underline{Step 3.}
Set $w(x,t):=\gamma W(x,t)+C_0 t$. 
By \eqref{eq:4.4}, \eqref{eq:5.2} and \eqref{eq:5.11} 
we have 
\begin{equation*}
\begin{split}
\Psi[w](x,t)
 & \le C\int_0^t \int_{{\bf R}^N} G(x-y,t-s) 
	\left( U(y,s)^p + \gamma^p W(y,s)^p + C_0^p s^p \right) \,dy\,ds\\
 & \le C\gamma^p W+CC_0^p t^{p+1}
\end{split}
\end{equation*}
for $(x,t)\in{\bf R}^N\times(0,\infty)$. 
Then we have: 
\begin{itemize}
  \item[(1)] 
  Let $\gamma>0$ and $T>0$ be sufficiently small. 
  Then  $\Psi[w](x,t)\le w$ in ${\bf R}^N\times[0,T)$, 
  that is, $w$ is a supersolution to problem~\eqref{eq:4.2} in ${\bf R}^N\times[0,T)$;
  \item[(2)] 
  Let $C_0=0$ and $\gamma>0$ be sufficiently small. 
  Then $\Psi[w](x,t)\le w$ in ${\bf R}^N\times[0,\infty)$, 
  that is, $w$ is a supersolution to problem~\eqref{eq:4.2} in ${\bf R}^N\times[0,\infty)$.
\end{itemize}
Thus the proof of Proposition~\ref{Proposition:5.1} is complete. 
$\Box$
\vspace{3pt}
\newline
Theorem~\ref{Theorem:1.4} with $p>p_*$ follows from 
Lemma~\ref{Lemma:4.1} and Proposition~\ref{Proposition:5.1}.  
\section{Proof of Theorem~\ref{Theorem:1.4} with $\boldsymbol{p=p_*}$}
Let $p=p_*$ and $\beta>0$. Let $A>1$ be as in Lemma~\ref{Lemma:4.3}. 
Set 
\begin{equation}
\label{eq:6.1}
W_*(x,t):= \int_0^t \int_{{\bf R}^N} G(x-y,t-s) g(t-s) h(f(y)) \,dy\,ds,
\end{equation}
where
\begin{equation}
\label{eq:6.2}
\begin{split}
 & f(x):= |x|^{-N} \left[\log\biggr(A+\frac{1}{|x|}\biggr)\right]^{-\frac{N}{\theta}}, \qquad 
g(\tau):= \left[\log\biggr(A+\frac{1}{\tau}\biggr) \right]^{-\beta},\\
 & h(\eta):= \eta \left[\log(A+\eta) \right]^\beta. 
\end{split}
\end{equation}
Then 
\begin{equation}
\label{eq:6.3}
C^{-1}|x|^{-N}\left[ \log\biggr(A+\frac{1}{|x|}\biggr) \right]^{-\frac{N}{\theta}+\beta} 
\leq h(f(x))\leq 
C|x|^{-N} \left[ \log\biggr(A+\frac{1}{|x|}\biggr) \right]^{-\frac{N}{\theta}+\beta}
\end{equation}
for $x\in{\bf R}^N$. 
Similarly to Section~5, 
Theorem~\ref{Theorem:1.4} with $p=p_*$ follows from Lemma~\ref{Lemma:4.1} and 
the following proposition. 
\begin{proposition}\label{Proposition:6.1}
Assume the same conditions as in Theorem~{\rm\ref{Theorem:1.4}}.
Let $p=p_*$ and $W_*$ be as in \eqref{eq:6.1}. 
Then there exist $\gamma_*=\gamma_*(N,\theta)>0$ and $\beta=\beta(N,\theta)>0$ 
such that, for any $C_0\geq0$, 
the function $w_*$ defined by 
$$
w_*(x,t):=\gamma_*W_*(x,t)+C_0t
$$ 
is a supersolution to problem~\eqref{eq:4.2} in ${\bf R}^N\times[0,T)$ for some $T>0$. 
\end{proposition}
{\bf Proof.}
Let $p=p_*$ and assume \eqref{eq:1.7}. 
Let $x\in{\bf R}^N\setminus\{0\}$, $t>0$ and $U$ be as in \eqref{eq:4.1}. 
We use the same notation as in Section~5. 
The proof is divided into 3 steps. 
\vspace{3pt}
\newline
\underline{Step 1.} We show 
\begin{equation}
\label{eq:6.4}
V[U^p](x,t)\le C\gamma^p[\log(A+t)]^{p-1}W_*(x,t)+CC_0^pt^{p+1}.
\end{equation}
It follows from \eqref{eq:1.7}, \eqref{eq:3.5} and \eqref{eq:4.1} that 
\begin{equation}
\label{eq:6.5}
U(x,t)\le C\gamma V[f](x,t)+C_0t. 
\end{equation}
Let $1<k<p$. 
By the H\"older inequality, \eqref{eq:6.2} and \eqref{eq:6.3} we have 
\begin{equation*}
\begin{split}
V[f](x,t) 
& \le V[h(f)](x,t)^{\frac{1}{k}}V[f^{\frac{k}{k-1}}h(f)^{-\frac{1}{k-1}}](x,t)^{1-\frac{1}{k}}\\
 & \le CV[h(f)](x,t)^{\frac{1}{p}}\\
 & \qquad\times\left(\int_0^t\int_{{\bf R}^N}G(x-y,t-s)|y|^{-N} 
\left[\log\biggr(A+\frac{1}{|y|}\biggr) \right]^{-\frac{N}{\theta}+\beta}\,dy\,ds\right)^{\frac{1}{k}-\frac{1}{p}}\\
 & \qquad\times\left(\int_0^t G(x-y,t-s) |y|^{-N} 
\left[\log\biggr(A+\frac{1}{|y|}\biggr) \right]^{-\frac{N}{\theta}-\frac{1}{k-1}\beta} \,dy\,ds\right)^{1-\frac{1}{k}}.
\end{split}
\end{equation*}
This together with Lemma~\ref{Lemma:4.3} and $p=p_*=N/(N-\theta)$ that 
\begin{equation}
\label{eq:6.6}
\begin{split}
V[f](x,t)
 & \le CV[h(f)](x,t)^{\frac{1}{p}}\\
 & \qquad\times\left(|x|^{\theta-N} 
	\left[\log\biggr(A+\frac{1}{|x|}\biggr) \right]^{-\frac{N}{\theta}+\beta+1} \log(A+t)\right)^{\frac{1}{k}-\frac{1}{p}}\\
 & \qquad\times\left(|x|^{\theta-N} 
	\left[\log\biggr(A+\frac{1}{|x|}\biggr) \right]^{-\frac{N}{\theta}-\frac{1}{k-1}\beta+1} 
	\log(A+t)\right)^{1-\frac{1}{k}}\\
 & \le C[\log(A+t)]^{1-\frac{1}{p}}|x|^{-\frac{\theta}{p}}\left[\log\biggr(A+\frac{1}{|x|}\biggr) \right]^{-\frac{\beta+1}{p}}
 V[h(f)](x,t)^{\frac{1}{p}}.
\end{split}
\end{equation}
By \eqref{eq:6.5} and \eqref{eq:6.6} we obtain 
$$
U(x,t)^p
\le C\gamma^p[\log(A+t)]^{p-1}\psi(x)
V[h(f)](x,t)+CC_0^pt^p,
$$
where
$$
\psi(x):=|x|^{-\theta}\left[\log\biggr(A+\frac{1}{|x|}\biggr) \right]^{-\beta-1}.
$$
This implies that 
\begin{equation}
\label{eq:6.7}
V[U^p](x,t) 
\le C\gamma^p[\log(A+t)]^{p-1}\tilde{I}(x,t)+CC_0^pt^{p+1},
\end{equation}
where 
\begin{equation}
\label{eq:6.8}
\begin{split}
 & \tilde{I}(x,t)\\
 & :=\int_0^t \int_{{\bf R}^N} h(f(z))
\left(\int_\tau^t \int_{{\bf R}^N} G(x-y,t-s) G(y-z,s-\tau)\psi(y)\,dy\,ds\right)\,dz\,d\tau\\
 & \,=\int_0^t \int_{{\bf R}^N}  G(x-z,t-\tau)h(f(z))\tilde{J}(x,z,t,\tau)\,dz\,d\tau
\end{split}
\end{equation} 
and 
$$
\tilde{J}(x,z,t,\tau):=
 \int_\tau^t \int_{{\bf R}^N}\frac{G(x-y,t-s) G(y-z,s-\tau)}{G(x-z,t-\tau)}\psi(y)\,dy\,ds.
$$
Similarly to \eqref{eq:5.7}, 
in the case of $\theta=2$, by Lemma~\ref{Lemma:4.5} we have 
\begin{equation}
\label{eq:6.9}
\begin{split}
\tilde{J}(x,z,t,\tau)
 & \le C\int_\tau^t \left\|\int_{{\bf R}^N}G(\cdot-y,\xi(s:t,\tau))\psi(y)\,dy\right\|_{L^\infty({\bf R}^N)}\,ds\\
 & \le C\int_\tau^t \xi(s:t,\tau)^{-1}\left[\log\biggr(A+\frac{1}{\xi(s:t,\tau)}\biggr)\right]^{-\beta-1}\,ds,
\end{split}
\end{equation}
where $\xi(s:t,\tau)$ is as in \eqref{eq:5.8}.
Since 
\begin{equation*}
\begin{split}
 & \frac{1}{2}(t-s)\le\xi(s:t,\tau)\le t-s\quad\,\,\mbox{for}\quad \frac{t+\tau}{2}\le s\le t,\\
 & \frac{1}{2}(s-\tau)\le\xi(s:t,\tau)\le s-\tau\quad\mbox{for}\quad \tau\le s\le\frac{t+\tau}{2},
\end{split}
\end{equation*}
by \eqref{eq:6.9} we obtain 
\begin{equation}
\label{eq:6.10}
\begin{split}
\tilde{J}(x,z,t,\tau)
 & \le C\int_{(\tau+t )/2}^t (t-s)^{-1}\left[\log\biggr(A+\frac{1}{t-s}\biggr)\right]^{-\beta-1}\,ds\\
 & \qquad+C\int_\tau^{(\tau+t )/2}(s-\tau)^{-1}\left[\log\biggr(A+\frac{1}{s-\tau}\biggr)\right]^{-\beta-1}\,ds\\
 & \le C\left[\log\biggr(A+\frac{1}{t-\tau}\biggr)\right]^{-\beta}.
\end{split}
\end{equation}
In the case of $0<\theta<2$, similarly to \eqref{eq:5.9}, 
it follows that
\begin{equation*}
\begin{split}
 & \tilde{J}(x,z,t,\tau)\\
 & \le C\int_\tau^t \left\|\int_{{\bf R}^N}G(\cdot-y,t-s)\psi(y)\,dy\right\|_{L^\infty({\bf R}^N)}\,ds\\
 & \qquad\qquad
 +C\int_\tau^t \left\|\int_{{\bf R}^N}G(\cdot-y,s-\tau)\psi(y)\,dy\right\|_{L^\infty({\bf R}^N)}\,ds\\
 & \le C\int_\tau^t (t-s)^{-1}\left[\log\biggr(A+\frac{1}{t-s}\biggr)\right]^{-\beta-1}\,ds
 +C\int_\tau^t (s-\tau)^{-1}\left[\log\biggr(A+\frac{1}{s-\tau}\biggr)\right]^{-\beta-1}\,ds\\
 & \le C\left[\log\biggr(A+\frac{1}{t-\tau}\biggr)\right]^{-\beta}.
\end{split}
\end{equation*}
This together with \eqref{eq:6.8} and \eqref{eq:6.10} implies that  
\begin{equation}
\label{eq:6.11}
\tilde{I}(x,t)\le C\int_0^t \int_{{\bf R}^N}  G(x-z,t-\tau)\left[\log\biggr(A+\frac{1}{t-\tau}\biggr)\right]^{-\beta} 
h(f(z))\,dz\,d\tau
=CW_*(x,t). 
\end{equation}
Then \eqref{eq:6.4} follows from \eqref{eq:6.7} and \eqref{eq:6.11}.
\vspace{3pt}
\newline
\underline{Step 2.} 
We prove
\begin{equation}\label{eq:6.12}
V[(W^*)^p](x,t) \leq C[\log(A+t)]^{p-1} W_*(x,t).
\end{equation}
It follows from Lemma \ref{Lemma:4.3} and $p=p_*=N/(N-\theta)$ that
\[
\begin{aligned}
	&\int_0^t \int_{{\bf R}^N} G(x-y,t-s) g(t-s)^\frac{p}{p-1} h(f(y)) \,dy\,ds \\
	&\leq 
	C\int_0^t \int_{{\bf R}^N} G(x-y,t-s)  
	\left[ \log\biggr(A+\frac{1}{t-\tau}\biggr) \right]^{-\frac{p}{p-1}\beta} 
	|y|^{-N}\left[ \log\biggr(A+\frac{1}{|y|}\biggr) \right]^{-\frac{N}{\theta}+\beta}  \,dy\,ds \\
	&\leq 
	C|x|^{\theta-N} 
	\left[\log\biggr(A+\frac{1}{|x|}\biggr) \right]^{-\frac{1}{p-1}\beta-\frac{N-\theta}{\theta} } 
	\log(A+t).
\end{aligned}
\]
Then it follows from the H\"older inequality and $p=p_*=N/(N-\theta)$ that 
\begin{equation*}
\begin{split}
 & W_*(x,t)^p
 \le \left(\int_0^t \int_{{\bf R}^N} G(x-y,t-s) g(t-s)^\frac{p}{p-1} h(f(y)) \,dy\,ds\right)^{p-1}V[h(f)](x,t)\\
 & \quad
\le C|x|^{-\theta} \left[\log\biggr(A+\frac{1}{|x|}\biggr) \right]^{-1-\beta} 
[\log(A+t)]^{p-1}\int_0^t \int_{{\bf R}^N} G(x-y,t-s) h(f(y)) \,dy\,ds.
\end{split}
\end{equation*}
By \eqref{eq:6.8} and \eqref{eq:6.11} we obtain 
$$
\int_0^t \int_{{\bf R}^N} G(x-y,t-s) W_*(y,s)^p \,dy\,ds
\le C[\log(A+t)]^{p-1}\tilde{I}(x,t)
\le C[\log(A+t)]^{p-1}W_*(x,t). 
$$
This implies \eqref{eq:6.12}. 
\vspace{3pt}
\newline
\underline{Step 3.}  
Set $w_*:=\gamma W_*+C_0t$. 
Let $\gamma>0$ and $T>0$ be sufficiently small. 
Then, by a similar argument to that of Step~3 in the proof of Theorem~\ref{Theorem:1.4} with $p>p_*$ 
we see that $w_*$ is a supersolution to problem~\eqref{eq:4.2} in ${\bf R}^N\times[0,T)$. 
Thus Proposition~\ref{Proposition:6.1} follows.
$\Box$
\vspace{5pt}

\noindent
{\bf Acknowledgments.}
The authors of this paper thank Professor Tatsuki Kawakami 
for informing them of Lemma~\ref{Lemma:4.4} in the case of $0<\theta<2$. 
The second author was supported in part 
by the Grant-in-Aid for Scientific Research (S) (No.~19H05599)
from Japan Society for the Promotion of Science. 
The third author was supported in part 
by the Early-Career Scientists  (No.~19K14567) 
from Japan Society for the Promotion of Science. 


\begin{thebibliography}{99}
\bibitem{A} 
R.~A.~Adams, 
{\it Sobolev Spaces}, 
Pure and Applied Mathematics, {\bf 65},  Academic Press, 1975. 

\bibitem{BLZ00}
C.~Bandle, H.~A.~Levine and Qi S.~Zhang,
Critical exponents of Fujita type for inhomogeneous parabolic equations and systems, 
J. Math. Anal. Appl. \textbf{251} (2000), 624--648. 

\bibitem{BP85} 
P.~Baras and M.~Pierre, 
Crit\`ere d'existence de solutions positives pour des \'equations semi-lin\'eaires non monotones, 
Ann. Inst. H. Poincar\'e Anal. Non Lin\'eaire {\bf 2} (1985), 185--212.

\bibitem{B}
G.~Bernard, 
Existence theorems for certain elliptic and parabolic semilinear equations, 
J. Math. Anal. Appl. {\bf 210} (1997), 755--776. 

\bibitem{BJ07}
K.~Bogdan and T.~Jakubowski, 
Estimates of heat kernel of fractional Laplacian perturbed by gradient operators, 
Comm. Math. Phys. \textbf{271} (2007), 179--198.

\bibitem{BSS03}
K.~Bogdan, A.~St\'os and P.~Sztonyk, 
Harnack inequality for stable processes on $d$-sets, 
Studia Math. {\bf 158}  (2003), 163--198.

\bibitem{Bu16}
C.~Bucur, 
Some observations on the Green function for the ball in the fractional Laplace framework, 
Commun. Pure Appl. Anal. \textbf{15} (2016), 657--699.

\bibitem{FK12}
A.~Z.~Fino and M.~Kirane, 
Qualitative properties of solutions to a time-space fractional evolution equation, 
Quart. Appl. Math. \textbf{70} (2012), 133--157. 

\bibitem{HI18}
K.~Hisa and K.~Ishige, 
Existence of solutions for a fractional semilinear parabolic equation with singular initial data,
Nonlinear Anal. \textbf{175} (2018), 108--132. 

\bibitem{HI19}
K.~Hisa and K.~Ishige, 
Solvability of the heat equation with a nonlinear boundary condition, 
SIAM J. Math. Anal. {\bf 51} (2019), 565--594.

\bibitem{IKO}
K.~Ishige, T.~Kawakami and S.~Okabe, 
Existence of solutions for a higher-order semilinear parabolic equation with singular initial data, 
arXiv:1909.05492v1. 


\bibitem{IS16}
K.~Ishige and R.~Sato, 
Heat equation with a nonlinear boundary condition and uniformly local $L^r$ spaces, 
Discrete Contin. Dyn. Syst. \textbf{36} (2016), 2627--2652. 

\bibitem{Ju05}
N.~Ju, 
The maximum principle and the global attractor 
for the dissipative 2D quasi-geostrophic equations,
Comm. Math. Phys. \textbf{255} (2005), 161--181.

\bibitem{KT16}
T.~Kan and J.~Takahashi, 
Time-dependent singularities in semilinear parabolic equations: behavior at the singularities, 
J. Differential Equations {\bf 260} (2016), 7278--7319. 

\bibitem{KT17}
T.~Kan and J.~Takahashi, 
Time-dependent singularities in semilinear parabolic equations: 
existence of solutions,
J. Differential Equations \textbf{263} (2017), 6384--6426.

\bibitem{KK04}
A.~G.~Kartsatos and V.~V.~Kurta, 
On blow-up results for solutions of inhomogeneous evolution equations and inequalities,
J. Math. Anal. Appl. \textbf{290} (2004), 76--85. 

\bibitem{KQ02}
M.~Kirane and M.~Qafsaoui, 
Global nonexistence for the Cauchy problem of some nonlinear reaction--diffusion systems,
J. Math. Anal. Appl. \textbf {268} (2002), 217--243.

\bibitem{L93}
T.-Y. Lee,
Some limit theorems for super-{B}rownian motion and semilinear differential equations. 
Ann. Probab. \textbf{21} (1993), 979--995.

\bibitem{MT06}
Y.~Maekawa and Y.~Terasawa, 
The Navier--Stokes equations with initial data in uniformly local $L^p$ spaces,
Differential Integral Equations \textbf{19} (2006), 369--400.

\bibitem{RS13}
J.~C. Robinson and M.~Sier\.{z}\polhk ega, 
Supersolutions for a class of semilinear heat equations,
Rev. Mat. Complut. \textbf{26} (2013), 341--360.

\bibitem{RS14}
X.~Ros-Oton and J.~Serra,
The Dirichlet problem for the fractional Laplacian: regularity up to the boundary,
J. Math. Pures Appl. \textbf{101} (2014), 275--302. 

\bibitem{Sugi75} 
S.~Sugitani, 
On nonexistence of global solutions for some nonlinear integral equations,
Osaka Math. J. {\bf 12} (1975), 45--51. 

\bibitem{Takahashi} 
J.~Takahashi, 
Solvability of a semilinear parabolic equation with measures as initial data, 
Geometric properties for parabolic and elliptic PDE's, 257--276, 
Springer Proc. Math. Stat., {\bf 176}, Springer, 2016.

\bibitem{We80}
F.~B.~Weissler, 
Local existence and nonexistence for semilinear parabolic equations in $L^p$, 
Indiana Univ. Math. J. {\bf 29} (1980), 79--102.

\bibitem{Zeng07}
X.~Z.~Zeng,
The critical exponents for the quasi-linear parabolic equations with inhomogeneous terms,
J. Math. Anal. Appl. \textbf{332} (2007), 1408--1424. 

\bibitem{Zhang98_1}
Qi S.~Zhang,
A new critical phenomenon for semilinear parabolic problems,
J. Math. Anal. Appl. \textbf{219} (1998), 125--139. 

\bibitem{Zhang98_2}
Qi S.~Zhang,
Blow up and global existence of solutions to an inhomogeneous parabolic system,
J. Differential Equations \textbf{147} (1998), 155--183. 

\bibitem{Zhang99}
Qi S.~Zhang,
Blow-up results for nonlinear parabolic equations on manifolds, 
Duke Math. J. \textbf{97} (1999), 515--539. 

\end{thebibliography}
\end{document}